\documentclass{amsart}

\usepackage{mathptmx}
\usepackage{amsthm}
\usepackage[ruled,vlined,linesnumbered]{algorithm2e}
\usepackage[T1]{fontenc}
\usepackage{booktabs}
\usepackage{graphicx}
\usepackage{dcolumn}
\usepackage{caption}
\usepackage{url}

\newtheorem{proposition}{Proposition}

\title[Resonance structure of the periodically forced S\&S model]
{Resonance structure of a periodically forced delay differential equation model for the El Niño--Southern Oscillation}

\author{S. Bolduc-St-Aubin}
\address{Department of Mathematics, University of Auckland, Private Bag 92019, Auckland 1142, New Zealand}
\email{sbol082@aucklanduni.ac.nz}

\author{P. Subramanian}
\address{Department of Mathematics, University of Auckland, Private Bag 92019, Auckland 1142, New Zealand}

\author{B. Krauskopf}
\address{Department of Mathematics, University of Auckland, Private Bag 92019, Auckland 1142, New Zealand}

\date{\today}

\begin{document}

\begin{abstract}
We study resonance phenomena in the periodically forced Suarez--Schopf delay differential equation, which is a conceptual climate model for the El Niño--Southern Oscillation (ENSO). The system serves as a prototypical forced delayed-action oscillator whose self-sustained oscillations, when subjected to periodic forcing, give rise to attracting invariant tori. We provide a comprehensive bifurcation analysis of both the unforced and the forced model; for the latter, we propose a method to compute the rotation number of normally hyperbolic attracting invariant tori. With it we show that resonance tongues in parameter space are organized by critical points of the graph of the rotation number, both along torus bifurcation curves and within the region of invariant tori. We also show that the resonance structure repeats for large delays, which constitutes a reappearance mechanism not previously reported in the literature. Furthermore, depending on the feedback strength, we find bistability between period-one orbits and invariant tori. This regime involves non-classical bifurcation sequences, including `saddle-node' and `gluing' bifurcations of tori.
\end{abstract}

\maketitle

\section{Introduction}\label{sec:level1}

The El Niño-Southern Oscillation (ENSO) is a complex and significant climate phenomenon that influences weather patterns globally  \cite{mcphaden2006enso}. Characterized by variations in sea surface temperatures (SSTs) and atmospheric pressure across the equatorial Pacific Ocean, ENSO manifests itself primarily in two phases: El Niño, with systematically warmer-than-average SSTs, and La Niña, with cooler-than-average SSTs in the eastern Pacific Ocean.  These events have a profound impact on global climate, altering weather patterns, ocean conditions, and atmospheric pressure systems. Such changes can have widespread consequences for agriculture, water resources, and socioeconomic systems worldwide  \cite{NIZAMANI2025101276}. 

\begin{figure}[t!] 
  \centering
  \includegraphics[width=10.5cm]{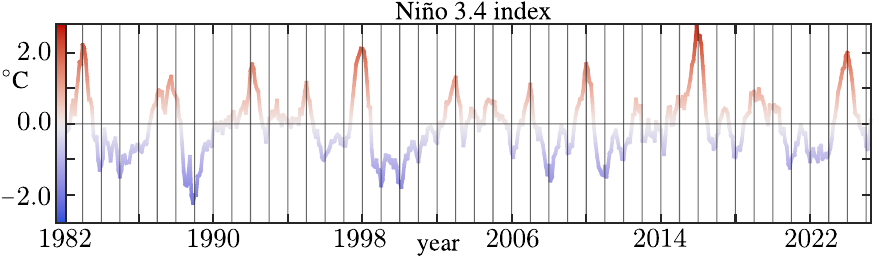}
  \caption{Niño 3.4 index time series taken from NOAA \cite{noaa_nino34_sst}, which records SSTs anomalies in the central equatorial Pacific (5°N–5°S and 170°W–120°W). El Niño is defined as the index exceeding +0.5°C for at least five consecutive overlapping 3-month periods; La Niña occurs when it falls below -0.5°C over the same duration  \cite{trenberth1997definition}.
  Peaks in the index occur near December, reflecting phase locking of ENSO variability to the seasonal cycle.
  }\label{fig:nino34}
\end{figure}

Figure~\ref{fig:nino34} shows the time series of the Ni\~no~3.4
index, a key indicator of El Niño and La Niña events. It illustrates the oscillatory nature of ENSO together with its marked
irregularity in both amplitude and period, but with a seasonal phase: El Niño and La Niña events typically peak around December  \cite{spencer2004role,stein2014enso}.

The earliest explanations of ENSO date back to Bjerknes~ \cite{bjerknes1966possible,bjerknes1969atmospheric}, who identified the \emph{Bjerknes positive feedback}: a warm SST anomaly in the eastern Pacific weakens the zonal SST gradient, reduces the Walker circulation~ \cite{Oliver2005} and easterly winds, and thereby promotes further eastward spreading of warm water. However, Bjerknes did not describe a mechanism capable of reversing this growth, making it clear that a \emph{negative feedback} must oppose the warming (or cooling) phase  \cite{wang2004understanding}.

Several conceptual frameworks formalize this interplay between positive and negative feedbacks, including the delayed action oscillator (DAO)  \cite{SuarezSchopf,battisti1989interannual}, the western Pacific oscillator  \cite{weisberg1997western,wang1999western}, the recharge–discharge oscillator  \cite{Jin1,Jin2}, and the advective–reflective oscillator  \cite{picaut1997advective}. Each emphasizes a different pathway through which the system restores balance  \cite{wang1999western}.

Conceptual frameworks like these condense the essential physics into manageable mathematical forms, facilitating analytical study and intuition building. This contrasts with large complex systems that model the oceans and atmosphere in considerable detail, and yet are much less amenable to mathematical discussion; see, e.g.,  \cite{DIJKSTRA2024133984} for a discussion of the role of conceptual climate models.

\begin{figure}
\includegraphics[width=10.9cm]{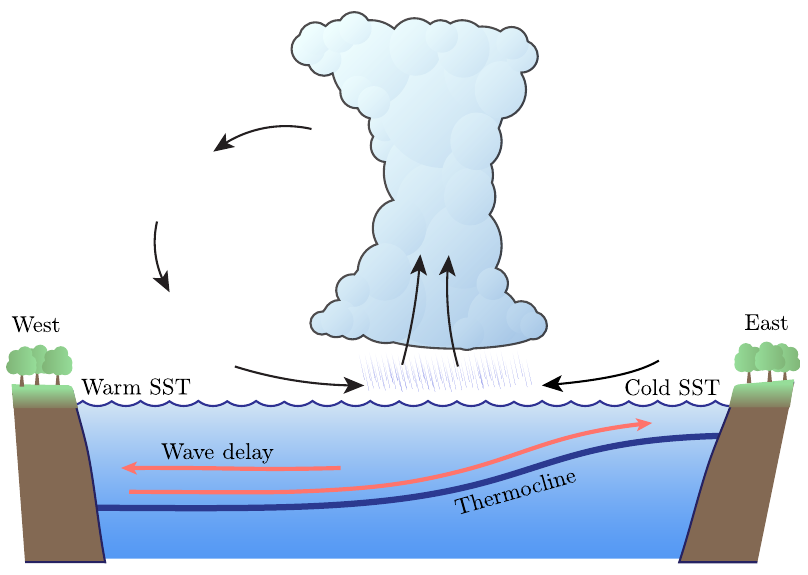}
\caption{\label{fig:feedback}  Sketch of a cross-section of the equatorial Pacific Ocean, where 
the blue line represents the thermocline. Oceanic waves (red arrows) travel west and east along the thermocline, creating a delayed negative feedback loop.} 
\end{figure}

In this work, we focus on the DAO paradigm introduced by Suarez and Schopf  \cite{SuarezSchopf}. Its underlying mechanism is illustrated schematically in Figure~\ref{fig:feedback} and works as follows.
Under normal conditions in the equatorial Pacific Ocean, easterly trade winds drive warm surface water toward the western equatorial Pacific, inducing upwelling in the east and establishing the Walker circulation. This zonal pattern produces warm, wet conditions in the west 
and cooler, drier conditions in the east. 
A perturbation that weakens the trade winds allows warm water to shift eastward.
This perturbation generates propagating oceanic waves that travel along the thermocline -- the transition zone between the warmer surface waters and the cooler deep waters.
Westward-propagating Rossby waves reflect at the western boundary and return as Kelvin waves, providing a \textit{negative delayed feedback} that 'terminates' the event. The perturbation also excites eastward-propagating Kelvin waves, providing a \textit{positive delayed feedback} that is not included in the original DAO formulation, but has been considered by other authors~ \cite{TSCJ}.

 Suarez and Schopf proposed a scalar delay differential equation for the ENSO DAO mechanism. In our notation, it takes the form
\begin{equation}
    \label{eq:ss}
    u'(t) = u(t) - u^3(t) - \alpha u(t-\tau),
\end{equation}
which we refer to as the \textit{Suarez and Schopf}, or simply, \textit{S\&S model} \cite{falkena2019derivation}. Here, $ u(t) $ denotes a non-dimensional SST anomaly in the equatorial eastern Pacific; the nondelayed linear term  captures the self-reinforcing nature of SST anomalies; and the cubic term acts as a nonlinear saturation mechanism. The delayed term $ -\alpha u(t - \tau) $ models the negative feedback, of strength $ \alpha $ and associated time delay $\tau$, representing the overall travel time of the Rossby and Kelvin waves along the thermocline across the Pacific basin.

We note that the cubic term  in the S\&S model~\eqref{eq:ss} was given little justification in the original work~ \cite{SuarezSchopf}. Later, Battisti and Hirst~ \cite{battisti1989interannual} provided qualitative support for this form, but without a mathematical derivation, making this model purely phenomenological. 
A significant advance was made by Falkena \emph{et al.}~ \cite{falkena2019derivation}, who derived S\&S--type equations from a spatially extended ENSO model via two complementary reduction procedures (a variation of constants and the Mori--Zwanzig formalism~ \cite{mori1965transport,zwanzig1973nonlinear}). Both approaches recover essentially the same form of the system, featuring a cubic nonlinearity, but the precise structure of the delayed feedback terms differs slightly. Importantly, the original S\&S model is retrieved from either formulation by setting a single parameter to zero, thereby placing it on a firm physical footing.

Motivated by the observation that ENSO events lock to the seasonal variation, several periodically forced DAO models have been considered and studied~ \cite{TSCJ,munnich1991study,Tony2026,boutle2007nino,ENSOjan2014}.
One of the most extensively studied conceptual ENSO models is the \emph{Ghil-Zaliapin-Thompson (GZT) model}~ \cite{GZT}, which combines periodic forcing and a delayed feedback mechanism with a switching function as the nonlinearity. Early investigations of the GZT model relied primarily on numerical
simulations~ \cite{GZT,GZT2}, while later work employed numerical
continuation techniques (\texttt{DDE-BifTool}~ \cite{ddebiftool}) to clarify the
structure of phase-locked solutions~ \cite{KKGZT,KKChenciner}.
Chaotic dynamics appear to be confined to relatively
small regions of parameter space where resonance tongues overlap~ \cite{KKGZT}.
Extensions of the model that incorporate state-dependent delay~ \cite{KKstatedep}
or asymmetric coupling function~ \cite{KKP_TSCJ} were shown to exhibit substantially larger regions of
resonance overlap and, correspondingly, more robust chaotic dynamics.
We note that analytical insight has been obtained from piecewise-defined or
reduced variants of the GZT model~ \cite{sgnmodel,Tony2026}.

Because the delayed negative feedback appears linearly, while nonlinearity acts locally through the cubic term, the S\&S model is structurally more transparent than the GZT model.
Moreover, in light of the recent derivations by Falkena
\emph{et al.}~ \cite{falkena2019derivation}, the S\&S model is no longer purely
phenomenological but can be interpreted as a systematically reduced
representation of spatially extended ENSO dynamics.
These considerations provide strong justification for revisiting the S\&S
framework to study the role of periodic forcing for observed dynamics. Namely, we study the system
\begin{equation}
\label{eq:pss}
    u'(t) = u(t) - u^3(t) - \alpha\,u(t-\tau) + c\,\cos(\omega t),
\end{equation}
where $\omega>0$ is the angular frequency and $c>0$ is the strength of the periodic forcing.
We refer to Eq.~\eqref{eq:pss} as the \textit{periodically forced S\&S} (pfS\&S) model. Its formulation combines a physically motivated delayed feedback structure with explicit seasonal forcing. 
This combination of terms makes Eq.~\eqref{eq:pss} relevant not only for climate modeling but also, more generally, for understanding the interaction between delayed feedback and periodic forcing.

In periodically forced delayed systems, such as the pfS\&S and GZT models, intrinsic oscillations associated with delayed feedback interact with the external seasonal forcing, allowing for motion on invariant tori. Trajectories on these tori can be characterized by the \emph{rotation number}, which represents the average number of rotations along one direction of the torus relative to the other. When the rotation number is rational, the dynamics is phase-locked; when the rotation number is irrational, the dynamics on the torus is quasiperiodic. In the former case, there are a stable and an unstable periodic orbit on the torus, and both form torus knots specified by the rational rotation number; see Section~\ref{section:compotationofrotation} for more details.
In a two-parameter bifurcation diagram, phase-locked orbits form resonance tongues whose interplay can lead to complicated bifurcation scenarios and chaotic behaviour~ \cite{chao,KKP_TSCJ}.

The pfS\&S  model~\eqref{eq:pss} has not yet been studied systematically. It was considered for angular forcing frequency $\omega=1$ by Anikushin and Romanov~ \cite{ANIKUSHIN2023133653}, who reported regimes with positive maximal Lyapunov exponents, indicating chaotic dynamics for selected parameter values. With computer-assisted proof techniques, Oishi and Sekine~ \cite{Oishi2021} studied Eq.~\eqref{eq:pss} for a wide range of angular frequencies $\omega$. They proved the existence of large, overlapping resonance tongues and strong resonances, suggesting rich dynamical scenarios.

\medskip

In this work we present a comprehensive bifurcation analysis of the pfS\&S model. In contrast to earlier work, we focus our attention on the case $\omega = 2\pi$, where the forcing period is one. The modeling idea here is that Eq.~\eqref{eq:pss} is expressed on a time scale of years. However, this is not obvious \emph{a priori}; because it depends on the time-rescaling factor used to obtained the non-dimensionalised Eq.~\eqref{eq:ss}; we come back to this question in Section~\ref{section:discussionandoutlook}.

We begin with a bifurcation study of the unforced S\&S model~\eqref{eq:ss} in the $(\tau,\alpha)$-plane. We perform a linear stability analysis of the equilibria, which yields analytic expressions for the pitchfork and Hopf bifurcations.
Bifurcations not captured by the linear stability analysis are studied with the continuation software \texttt{DDE-BifTool} \cite{ddebiftool}. Specifically, we find and continue saddle-node bifurcations of periodic orbits and homoclinic/gluing bifurcations. For sufficiently large delays, the only attracting invariant objects are periodic orbits.

We then analyse the bifurcations of the pfS\&S model~\eqref{eq:pss} in the $(\tau,c)$-plane for the representative values of $\alpha = 1.25$ and $\alpha = 0.75$. The interaction between self-sustained oscillations and external forcing gives rise to invariant tori, and a key tool to study these dynamics is the rotation number. To compute it, we adapt the method of Bailie~\emph{et al.}~ \cite{bailie2024resonancetongues,John} --- developed for periodically forced planar ODEs --- to the delay setting, where the infinite-dimensional phase space requires suitable projections; this is described in Section~\ref{section:compotationofrotation}. We show that resonance tongues in the $(\tau,c)$--plane are very thin, indicating weak locking. Specifically, no resonance tongues overlap nor are period-doubling cascades observed. Nevertheless, interesting resonance structures appear: resonance tongues are organized not only by extrema of the rotation-number along torus bifurcation curves~ \cite{terrien}, but also by critical points in the interior of the  region with attracting invariant tori; this is consistent with recent observations in a periodically forced planar system~ \cite{John}. For large delays, the resonance structure repeats periodically as described by rational map; at the level of the individual locked orbits, this corresponds to additional winding of the locked orbits along the longitudinal direction of the invariant torus. Moreover, symmetry-breaking bifurcations lead to complicated scenarios, including `folds of tori' and `gluing bifurcations of tori', in which invariant tori break up through a complex sequence of bifurcations. Despite this complexity, the overall dynamics — as studied via the headpoints on the stroboscopic section — closely parallel those of the unforced case.

We also briefly examine a periodically forced version of a physically motivated S\&S-type equation derived by Falkena \emph{et al}.~ \cite{falkena2019derivation} and show that it has effectively the same bifurcations and resonance structures; see Section~\ref{section:discussionandoutlook}. The only difference is that the resonance organization is richer: we find more critical points of the rotation number $\rho$ both along the torus bifurcation curves and within the interior of the region of invariant tori.

\medskip
The paper is organized as follows. 
Section~\ref{section:background} provides background on DDEs, dynamics on invariant tori, rotation numbers, and numerical methods for DDEs. 
In Section~\ref{section:II} we analyze the autonomous S\&S model~\eqref{eq:ss}. 
In Section~\ref{section:compotationofrotation}, we show that the pfS\&S model admits motions to attracting invariant tori and introduce an algorithm to compute their rotation numbers. 
Section~\ref{section:bif_pss} presents a bifurcation analysis of the pfS\&S model~\eqref{eq:pss} with angular frequency $\omega = 2\pi$ for the cases $\alpha = 1.25$ and $\alpha = 0.75$.
In Section~\ref{section:global}, we study the reappearance mechanisms of the resonance tongues for large $\tau$ values.
Finally, in Section~\ref{section:discussionandoutlook} we briefly present for comparison the resonance structure of the periodically forced version of a model by Falkena~\emph{et al.}~ \cite{falkena2019derivation}, discuss the role of varying the forcing frequency $\omega$, and outline future directions.

\section{Background}
\label{section:background}

We recall here some material underlying the analysis and numerical framework of our study.

\subsection{Delay Differential Equations Formulation}
\label{section:dde}
Model \eqref{eq:ss} is an example of autonomous \emph{delay differential equations} (DDEs), where the time derivative at the current time depends explicitly on the state of the system at earlier times. More specifically, Eq.~\eqref{eq:pss} is a scalar DDE with a single constant delay $\tau$ of the form
\begin{equation}   \label{eqn:DDE}
    {u}'(t)=f(u(t),u(t-\tau)),\quad f:\mathbb{R}\times\mathbb{R}\to\mathbb{R}.
    \end{equation}
An entire segment of the past trajectory is required as an initial condition, namely an initial function that specifies the state over $[-\tau,0]$. 
Consequently, the phase space of a DDE of the form \eqref{eqn:DDE} is the infinite-dimensional space $\mathcal{C}:=\mathcal{C}([-\tau,0],\mathbb{R})$ of continuous functions over $[-\tau,0]$ with values in $\mathbb{R}$.
Given the right-hand side of \eqref{eqn:DDE} is locally Lipschitz continuous,
the initial value problem with continuous initial data in $\mathcal{C}$
admits a unique continuously differentiable solution $u(t)$ for some $t \geq 0$; see, e.g.,  \cite{hale1993,smithdde,diekmann2012delay}.
For an initial function segment $\varphi\in\mathcal{C}$, the corresponding solution operator defines a continuous \emph{semiflow}
\[
S(t) : \mathcal{C} \to \mathcal{C}, \qquad S(t)(\varphi) = u_t,
\]
for some $t\geq 0$,
on the Banach space \(\mathcal{C}\) equipped with the supremum norm; here $u(t)=u^\varphi(t)$ denotes the unique solution of \eqref{eqn:DDE} with initial condition $\varphi$, and
$u_t\in\mathcal{C}$ is the state of the system at time $t$, defined in the usual way as 
\[u_t(s) := u(t+s),\qquad s\in[-\tau,0].\] In particular, $u_t(0)=u(t)$ and $u_t(-\tau)=u(t-\tau)$.

The infinite-dimensional phase space allows scalar DDEs with constant delays to exhibit a surprisingly rich range of dynamics, including sustained oscillations, chaotic behavior, and strong sensitivity to initial conditions; see, e.g., \cite{Duruisseaux_2022}.  
 At the same time, the local stability theory of equilibria and periodic orbits is entirely analogous to that of \emph{ordinary differential equations} (ODEs) since, as in ODEs, linearisations around equilibria and periodic orbits have at most finitely many unstable eigendirections \cite{diekmann2012delay}.
 Hence, classical bifurcation scenarios extend naturally to DDEs.

Their apparent simplicity in formulation paired with richness in behavior makes  DDEs particularly appealing for exploring complex dynamical phenomena without resorting to overly large or computationally demanding models \cite{chao}. 

\medskip
In the presence of periodic forcing,
the DDE \eqref{eqn:DDE} 
takes the form
\begin{equation}   \label{eqn:pDDE}
    u'(t)=f(u(t),u(t-\tau),t),\qquad 
    f:\mathbb{R}\times\mathbb{R}\times\mathbb{S}^1\to\mathbb{R},
\end{equation}
where $f(\cdot,\cdot,t+1)=f(\cdot,\cdot,t)$ is a periodic function whose period can be assumed to be $1$ without loss of generality.
Since system~\eqref{eqn:pDDE} is non-autonomous, the semiflow $S(t)$ 
no longer applies. Instead, one considers the two-parameter solution operator
\[
S(t,t_0)\colon\mathcal{C}\to\mathcal{C},
\qquad
S(t,t_0)(\varphi)=u_t,
\quad t\geq t_0,
\]
where $u(t)=u^{\varphi,t_0}(t)$ is the unique solution of~\eqref{eqn:pDDE} 
with $u_{t_0}=\varphi$; see~\cite{hale1993,diekmann2012delay}.
The phase space is accordingly extended to $\mathcal{C}\times\mathbb{S}^1$,
where $\mathbb{S}^1$ denotes the circle.

The theory of periodically forced DDEs is similar to that of periodically forced ODEs. 
In systems of the form \eqref{eqn:pDDE}, there are generally no equilibria, and the simplest attractors are period-one orbits. In particular, interactions between intrinsic oscillations and the periodic forcing generally leads to dynamics on a two-dimensional invariant torus $\mathcal{T}\cong\mathbb{S}^1\times\mathbb{S}^1$. As for ODEs, trajectories 
lie on $\mathcal{T}$. For the DDE~\eqref{eqn:pDDE}, this means that full history segments  lie on the two-dimensional invariant torus $\mathcal{T}$.

A practical way to study the invariant objects of system~\eqref{eqn:pDDE} is to exploit the periodic forcing. We define the 
\emph{stroboscopic map}
\begin{equation}
    \label{eq:stroboscopic}
    \mathcal{P}(\varphi;t_0)=S(t_0+1,t_0)(\varphi),
\end{equation}
which maps the history segment $\varphi=u_{t_0}$ at time $t_0$
to the segment $u_{t_0+1}$
one period later; the stroboscopic map acts on the infinite-dimensional space $\mathcal{C}$ by identifying $t_0$ and $t_0+1$ (since $f(\cdot,\cdot,t)$ in \eqref{eqn:pDDE} has period 1). 
Note that \eqref{eq:stroboscopic} is an example of a Poincar\'e map for the section $t=t_0$; in practice, one normally chooses $t_0=0$, and we simply write $\mathcal{P}(\varphi)=\mathcal{P}(\varphi;0)$ from now on. Period-one orbits of system~\eqref{eqn:pDDE} correspond to fixed points of $\mathcal{P}$ and, more generally, period-$n$ orbits of \eqref{eqn:pDDE} correspond to period-$n$ points of $\mathcal{P}$, that is, fixed points of $\mathcal{P}^n$. Because the map $\mathcal{P}$ acts on the infinite dimensional function space $\mathcal{C}$, projection methods are generally needed to visualize and study its invariant objects; see Section~\ref{section:compotationofrotation}.

\subsection{Motions on invariant tori and the rotation number}
\label{section:background_rot}
The stroboscopic map $\mathcal{P}$ provides a natural framework 
for characterizing the dynamics on an invariant torus $\mathcal{T}$ via the \emph{rotation number}. Note that $\mathcal{T}$ is two-dimensional, even though it lies in the infinite-dimensional phase space $\mathcal{C}\times\mathbb{S}^1$.

For periodically forced planar ODEs, an invariant torus $\mathcal{T}$ defines an invariant closed curve $\mathcal{I}$ of the stroboscopic map $\mathcal{P}$, on which $\mathcal{P}$ acts as a circle map. 
The \emph{rotation number} $\rho$ is 
the limit of the average angular displacement per iterate along the 
invariant closed curve $\mathcal{I}$; 
see, e.g., \cite{Kuznetsov}. 
The rotation number can equivalently be determined as the average of the oriented angular increments between consecutive iterates $z_k$ of $\mathcal{P}$ on $\mathcal{I}$:
\begin{equation}
    \label{eq:rotpractice}
    \rho = \lim_{N\to\infty} \frac{1}{2\pi N}\sum_{k=0}^{N-1}\angle(z_k,{Z}, z_{k+1}),
\end{equation}
where $\angle(z_k,{Z}, z_{k+1})$ denotes the oriented angle from $z_k$ to $z_{k+1}$, taken in $(-\pi,\pi]$, with respect to a suitable reference point ${Z}$ in the interior of $\mathcal{I}$; note that $\rho$
 is independent of the choice of reference point, provided ${Z}$ lies in the interior of (a projection of) the invariant closed curve $\mathcal{I}$; see \cite{John}.

\smallskip

The rotation number $\rho$ distinguishes two qualitatively different regimes:
\begin{itemize}
    \item \textbf{locked dynamics}: when $\rho =p/q\in\mathbb{Q}$ (with $p$ and $q$ relative prime) then there are an attracting and a repelling periodic orbit on the torus, each of which forms a $p$:$q$ torus knot on~$\mathcal{T}$.  
    \item \textbf{quasiperiodic dynamics}: when $\rho\notin\mathbb{Q}$ then every trajectory on $\mathcal T$ is dense, that is, never closes up.
\end{itemize}

Locked solutions are organized into resonance (or Arnold) tongues in parameter space, which are bounded by saddle-node bifurcations of periodic orbits~ \cite{Kuznetsov}. 
Specifically, a $p$:$q$ resonance tongue emanates from a point on a torus bifurcation curve \textsf{T} where the critical Floquet multipliers are $e^{\pm i2\pi p/q}$. 
The critical Floquet multipliers change continuously along \textsf{T}; moreover, resonance tongues are organized locally by the Farey tree~ \cite{KKGZT}: the largest resonance tongue between those with rotation numbers $p/q$ and $\ell/m$ is the one with rotation number $(p+\ell)/(q+m)$. Near the torus bifurcation curve \textsf{T}, the resonance tongues do not overlap; however, further away from \textsf{T}, neighboring resonance tongues may overlap, with associated breakdown of tori and possibly chaotic dynamics.

For DDEs, the invariant torus remains two-dimensional despite the fact that it now lies in an infinite-dimensional phase space. This is why the theory of resonance tongues carries over. Indeed, dynamics on invariant tori have been studied extensively for certain DDE models; see, e.g., \cite{KKGZT,Duruisseaux_2022,Calleja2017ResonancePI,terrien,ruschel2025strongly}. This was achieved by finding the bounding curves \textsf{SN} of saddle-node bifurcation of periodic orbits of selected resonance tongues with continuation tools.
The rotation number $\rho$, on the other hand, has not been considered for DDEs because computing it requires additional care, as we will see in Section~\ref{section:compotationofrotation}.

\subsection{Numerical tools}
Several numerical tools are available for DDEs, and we use here a combination of numerical integration, continuation, and the computation of the rotation number $\rho$.

Time integration is important and has several uses in this work. We make use of the \texttt{dde23}~ \cite{dde23} solver in \texttt{MATLAB} to compute the unstable manifold of an equilibrium, where our trajectories start in the unstable eigendirection of the equilibrium (which is provided by \texttt{DDE-BifTool}). We also compute the unstable manifold of a periodic orbit in a Poincar\'{e} section in a similar fashion; specifically, with the method detailed in \cite{krauskopf2003computing}, which requires several starting trajectories.
Furthermore, we use \texttt{dde23} for finding stable solutions on tori, which are then used to find and continue the saddle-node bifurcations that bound resonance tongues in parameter space. To compute the rotation number, as detailed in 
Section~\ref{section:compotationofrotation}, we perform extensive iterations by numerical integration of the stroboscopic map $\mathcal{P}$ over a large parameter grid. For efficiency, we employ an Adams--Bashforth method~ \cite{durran1991third} implemented in \textsc{Julia}~ \cite{bezanson2017julia}.

\sloppy{An essential tool we employ is the continuation software \texttt{DDE-BifTool}~ \cite{ddebiftool} in \texttt{MATLAB}~ \cite{MATLAB2024b}. The package \texttt{DDE-BifTool} is designed for autonomous systems.
To apply it also to the non-autonomous system~\eqref{eqn:pDDE}, such as system \eqref{eq:pss}, we rewrite it as an autonomous system by introducing auxiliary equations that generate the periodic forcing. This is achieved by increasing the dimension of the system by two, to implement the Hopf normal form with an attracting unit circle; see \cite{KKGZT}. \texttt{DDE-BifTool} enables the numerical continuation of equilibria and periodic orbits, as well as the continuation of bifurcations in one and two parameters. It also implements normal form computations, allowing, for instance, the determination of the criticality of Hopf bifurcations. The linear stability analysis also provides eigenvalues and associated eigenfunctions of the linearized problem. Furthermore, we continue saddle-node bifurcation of periodic orbits bounding resonance tongues and locate their root points along torus bifurcation curves by tracking the critical Floquet multipliers.}

\section{The Autonomous S\&S model}
\label{section:II}

The S\&S model \eqref{eq:ss} is dissipative \cite{Anikushin2022}. In particular, for $\alpha\in(0,1)$, all trajectories of Eq.~\eqref{eq:ss} are eventually attracted to the set $\{\varphi\in\mathcal{C}:\|\varphi\|_\infty\leq\sqrt{1+\alpha}\}$; see \cite{ANIKUSHIN2023133653}.
Moreover, Eq.~\eqref{eq:ss} has a natural  \(\mathbb{Z}_2\)-symmetry: if \(u(t)\) is a solution then so is \(-u(t)\).  Consequently, the global attractor (which exists due to dissipativity) is invariant under the involution \(u\mapsto -u\), or it consists of two parts that are each others image under the symmetry; the same holds for  any other invariant sets.

\subsection{Linear stability analysis}
\label{subsection:linearAnalysis}

Equilibria of DDEs are constant functions in the space $\mathcal{C}$. Their stability analysis follows the classical approach for ODEs of linearizing the equation around the equilibrium. However, due to the presence of a delay term, the characteriztic equation is transcendental rather than polynomial~ \cite{smithdde}.
Nevertheless, bifurcations occur as parameters vary when characteriztic roots cross the imaginary axis. 

\begin{proposition}
\label{prop:hopf}
[Local bifurcations] 

\noindent
For the autonomous S\&S model~\eqref{eq:ss} with parameters
$\tau,\alpha>0$, we have the following.
\begin{enumerate}
    \item[\textnormal{(i)}]
    For $\alpha>1$, 
    \[p_0\equiv0\]
    is the only equilibrium, which is attracting for small $\tau$.
    \item[\textnormal{(ii)}]
    The equilibrium $p_0$ undergoes a pitchfork bifurcation $\textsf{P}(\tau,\alpha)$ along the line given by $\alpha = 1$; as $\alpha$ decreases through $1$, the equilibria
    \[p^\pm = \pm\sqrt{1-\alpha}\]
    emerge.
    
    \item[\textnormal{(iii)}]
    For $0<\alpha<1$, the equilibrium $p_0$ is a saddle  with a one-dimensional unstable manifold $W^u(p_0)$. The equilibria $p^\pm$ are attracting for small $\tau$.
    
    \item[\textnormal{(iv)}]
   The equilibrium $p_0$ undergoes Hopf bifurcations $\textsf{H}_k(\tau_k,\alpha)$ at
    \begin{equation}
    \label{eq:tauk-statement}
        \tau_k
        =
        \frac{
            \arccos(1/\alpha)+2\pi k
        }
        {\sqrt{\alpha^2-1}},
        \qquad k\in\mathbb{N}_0.
    \end{equation}
    Any periodic orbit $u(t)$ born at these Hopf bifurcations is symmetric, and satisfies
    \begin{equation}
        \label{eqn:symmetry}
        u(t) = -u(t+\tfrac{T}{2}),
    \end{equation}
    where $T$ is its (minimal) period. The principal branch $\textsf{H}_0(\tau_0,\alpha)$ marks the first loss of stability of $p_0$ as $\tau$ increases; $\textsf{H}_0(\tau_0,\alpha)$ is supercritical and the bifurcating periodic orbits is attracting for small $\tau$. The curves $\textsf{H}_k(\tau_k,\alpha)$, $k\geq 1$, lie in the
    region $\alpha>1$ and have $\alpha=1$ as an asymptote as $\tau\to\infty$, and 
    $\tau=0$ as an asymptote as $\alpha\to\infty$.

    \item[\textnormal{(v)}]
    The equilibria $p^\pm$ undergo Hopf
    bifurcations $\textsf{H}_k^\pm(\tau_k^\pm,\alpha)$ at
    \begin{equation}
    \label{eq:taukdagger-statement}
        \tau_k^\pm
        =
        \frac{
            \arccos\!\big((3\alpha-2)/\alpha\big)
            +2\pi k
        }
        {2\sqrt{(1-2\alpha)(\alpha-1)}} ,
        \qquad k\in\mathbb{N}_0.
    \end{equation}
    All periodic orbits born at these Hopf bifurcations come in pairs, related to each other
    by the symmetry~\eqref{eqn:symmetry}. The principal branch $\textsf{H}_0^\pm(\tau_0^\pm,\alpha)$
    marks the first loss of stability of $p^\pm$ as $\tau$ increases; $\textsf{H}_0^\pm(\tau_0^\pm,\alpha)$ is subcritical and the bifurcating periodic orbit is a saddle. The curves $\textsf{H}_k^\pm(\tau_k^\pm,\alpha)$, $k\geq 1$, exist in the region $\tfrac{1}{2}<\alpha<1$,
    and have asymptotes $\alpha=\tfrac{1}{2}$ and $\alpha=1$ as $\tau\to\infty$.
    
    \item[\textnormal{(vi)}]
    There is a point \textsf{DZ} of double-zero bifurcation at $(\tau,\alpha) = (1,1)$ on the curve \textsf{P}, 
   here the linearisation has a double eigenvalue at zero, analogous to a Bogdanov--Takens bifurcation, but this occurs here in a system with $\mathbb{Z}_2$-symmetry; see e.g. \cite{sieber2004bifurcation}. The curves $\textsf{H}_0$, and $\textsf{H}_0^\pm$ emerge
    from the point \textsf{DZ}. Because of the $\mathbb{Z}_2$-symmetry of Eq.~\eqref{eq:ss}, the double-zero point \textsf{DZ} has the unfolding of the 1:2 resonance  with $s=-1$ in Chapter~8 of \cite{Kuznetsov}.
    
\end{enumerate}
\end{proposition}

\noindent 
For the proof see
Appendix~\ref{appendix:local}.

\subsection{One-parameter bifurcation diagrams}
\label{subsection:numericalanalysis}

\begin{figure}[htp] 
  \centering
  \includegraphics[width=9.09cm]
  {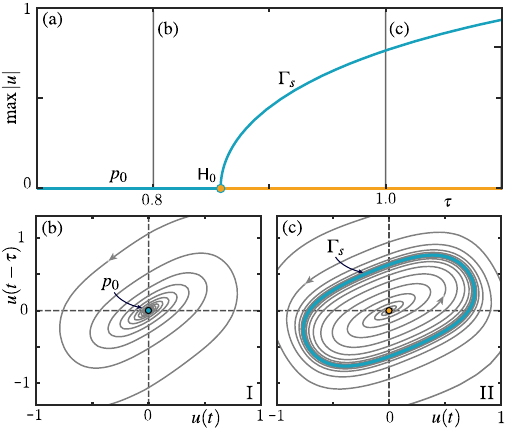}
  \caption{
    One-parameter bifurcation diagram in $\tau$ for $\alpha = 1.25$. Panel (a) shows the bifurcation diagram for $\tau\in[0.7,1.1]$ in terms of the max-norm; attracting objects are blue, and saddles are orange. The equilibrium $p_0$ loses stability at $\textsf{H}_0$, yielding the attracting periodic orbit $\Gamma_s$. Remaining panels show phase portraits in the $(u(t),u(t-\tau))$-plane at $\tau=0.8$ (b) and $\tau=1.0$ (c), with representative trajectories (grey).}\label{fig:oneparameteralpha1255}
\end{figure}

\begin{figure}[t!] 
  \centering
  \includegraphics[width=9cm]{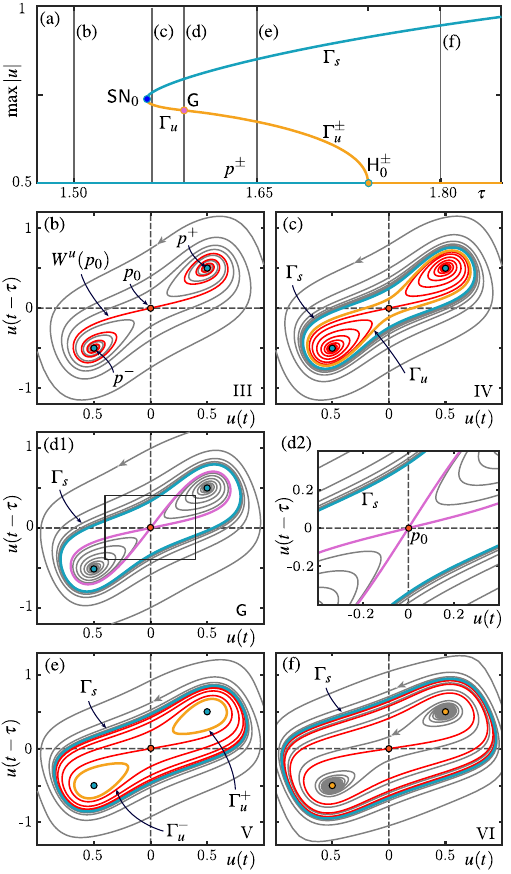}
  \caption{
One-parameter bifurcation diagram in $\tau$ for $\alpha = 0.75$. Panel (a) shows the bifurcation diagram for $\tau \in [1.475, 1.85]$ in terms of the max-norm, with saddle-node bifurcation $\textsf{SN}_0$, gluing bifurcation \textsf{G} and Hopf bifurcations $\textsf{H}_0^\pm$. Remaining panels show phase portraits in the $(u(t),u(t-\tau))$-plane at $\tau=1.5$ (b), $1.5656$ (c),  $1.5903$ (d1), $1.65$ (e) and $1.8$ (f). Panel~(d2) is an enlargement of panel~(d1). 
Shown are: equilibrium $p_0$ (orange); equilibria $p^\pm$ (blue when attracting, orange when saddle); attracting periodic orbit $\Gamma_s$ (blue); saddle periodic orbits $\Gamma_u$ and $\Gamma_u^\pm$ (orange); unstable manifold $W^u(p_0)$ (red); homoclinic connection (purple); and representative trajectories (grey), with arrows showing the direction of the semiflow.}\label{fig:oneparameteralpha075}
\end{figure}

To illustrate Proposition~\ref{prop:hopf}, Figure~\ref{fig:oneparameteralpha1255} and Figure~\ref{fig:oneparameteralpha075} show one-parameter bifurcation diagrams in 
$\tau$ for fixed values of $\alpha=1.25$ and $\alpha=0.75$, respectively. 
In each figure, panel~(a) presents the one-parameter bifurcation diagram in $\tau$ with branches of solutions and associated bifurcations. 
The remaining panels show phase portraits in projections onto the $(u(t),u(t-\tau))$-plane, with equilibria and periodic orbits, and some representative trajectories. For the case $\alpha=0.75$, we also show the one-dimensional unstable manifold $W^u(p_0)$. In this projection, the $\mathbb{Z}_2$-symmetry is evident: orbits are either invariant under a  rotation by $\pi$ about the origin, or they occur in symmetric pairs.

Figure~\ref{fig:oneparameteralpha1255}(a) shows the one-parameter bifurcation diagram in $\tau$ for $\alpha=1.25$, which has two regimes, labeled I and II. In regime I for small delay $\tau$, the equilibrium $p_0$ is a global attractor, as seen in panel~(b), where trajectories spiral toward $p_0$. The supercritical Hopf bifurcation $\textsf{H}_0$ occurs at $\tau_0\approx0.8580$. 
From this bifurcation point emerges a branch of symmetric attracting periodic orbits $\Gamma_s$. The existence of $\Gamma_s$ is a feature of regime II: as panel~(c) shows, all trajectories converge to this symmetric periodic orbit.

Figure~\ref{fig:oneparameteralpha075}(a) shows the bifurcation diagram in $\tau$ for $\alpha = 0.75$. Compared with the previous case, the dynamics here are more intricate, and we distinguish the five distinct regimes III--VI. 
Panel~(b) shows the phase portrait of regime III where the system exhibits the expected behavior following the pitchfork bifurcation~\textsf{P}: the symmetric equilibria $p^\pm$ are attracting, while $p_0$ is a saddle. The one-dimensional unstable manifold $W^u(p_0)$ has two sides, and one converges to $p^+$ and, by symmetry, the other converges to $p^-$. All other trajectories (not on the stable manifold of $p_0$)  converge to either $p^+$ or $p^-$.
As $\tau$ is increased, we encounter a bifurcation not covered by Proposition~\ref{prop:hopf}: a saddle-node bifurcation of periodic orbits $\textsf{SN}_0$. The phase portrait in the ensuing regime IV, shown in panel~(c), depicts an attracting outer orbit $\Gamma_s$ and a repelling inner orbit $\Gamma_u$. Any trajectory inside $\Gamma_u$ converges to either $p^+$ or $p^-$, while trajectories outside $\Gamma_u$ converge to $\Gamma_s$. In particular, the unstable manifold $W^u(p_0)$  still converges to $p^+$ and $p^-$.
At $\tau\approx 1.5903$, the inner periodic orbit $\Gamma_u$ has shrunk and now connects to the saddle equilibrium $p_0$, forming the figure-eight structure in Figure~\ref{fig:oneparameteralpha075}(d1). This homoclinic bifurcation is known as a \emph{gluing bifurcation}, which we refer to as \textsf{G}. Here, both sides of the unstable manifold $W^u(p_0)$ return to the equilibrium $p_0$; see the enlarged view in panel~(d2).
The result of this gluing bifurcation is regime~V, shown in panel~(e), where a pair of repelling periodic orbits $\Gamma_u^\pm$ exists. Both sides of $W^u(p_0)$ now converge to the attracting periodic orbit $\Gamma_s$, as do all other trajectories, except those inside the region enclosed by the pair $\Gamma_u^\pm$, which converge to either $p^+$ or $p^-$.
For increasing $\tau$, the repelling periodic orbits $\Gamma_u^\pm$ shrink and eventually disappear at the subcritical Hopf bifurcations $\textsf{H}_0^\pm$ at $\tau_0^\pm \approx 1.741$; its first Lyapunov coefficient is $l_1\approx 2.9574$, indicating that the Hopf bifurcation is indeed subcritical. Beyond this point, in regime~VI, shown in panel~(f), the only remaining attractor is the periodic orbit $\Gamma_s$.

\begin{figure}[t!] 
  \centering
  \includegraphics[width=9.09cm]{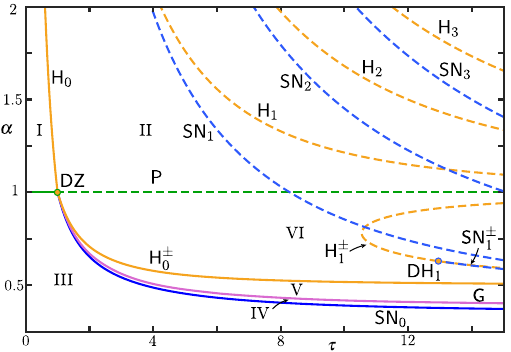}
  \caption{Bifurcation diagram of the S\&S model~\eqref{eq:ss} in the $(\tau,\alpha)$-plane showing the regimes I--VI; the curves of pitchfork bifurcation $\textsf{P}$, Hopf bifurcations $\textsf{H}_k$, $\textsf{H}_k^\pm$, fold bifurcations of periodic orbits $\textsf{SN}_k$, $\textsf{SN}_k^\pm$,  and homoclinic/gluing $\textsf{G}$; and the codimension-two points $\textsf{DZ}$ of double-zero bifurcation and $\textsf{DH}_1$ of degenerate Hopf bifurcation.
  }\label{fig:taualpha}
\end{figure}

\subsection{Two-parameter bifurcation diagram in the $(\tau,\alpha)$-plane}
\label{subsub:two}

Figure~\ref{fig:taualpha} presents the bifurcation diagram of system \eqref{eq:pss} in the $(\tau, \alpha)$-parameter plane, over the ranges $\tau \in [0,15]$ and $\alpha \in [0,2]$. Shown are the bifurcation curves \textsf{P}, $\textsf{H}_k$ and $\textsf{H}_k^\pm$ from Proposition~\ref{prop:hopf}, together with the additional curves \textsf{G} and $\textsf{SN}_0$ identified in Figure~\ref{fig:oneparameteralpha075}. Solid curves in Figure~\ref{fig:taualpha} indicate bifurcations involving at least one attracting invariant object; they separate the plane into the regions I--VI defined in Section~\ref{subsection:numericalanalysis}. Dashed curves represent bifurcations involving only saddle-type invariant objects; although the attractors themselves remain unchanged, saddle invariant objects may change or appear.

The pitchfork bifurcation \textsf{P} splits the $(\tau,\alpha)$-plane into two regions, with the equilibria $p^\pm$ existing only for $0<\alpha<1$. Consequently, regions I--II are located above the line $\alpha = 1$, whereas regions III--VI are below it.
On \textsf{P}, at $(\tau,\alpha)=(1,1)$, lies the double-zero point \textsf{DZ}, from which the curves $\textsf{H}_0$ and $\textsf{H}_0^\pm$ emanate, as stated in Proposition~\ref{prop:hopf}. The principal Hopf bifurcation curve $\textsf{H}_0$ separates regions I and II and approaches $\tau = 0$ as $\alpha \to \infty$.
The Hopf bifurcation curve $\textsf{H}_0^\pm$ tends to $\alpha = \tfrac12$ as $\tau \to \infty$. From \textsf{DZ} also emerge the saddle-node of periodic orbit curve $\textsf{SN}_0$ and the gluing bifurcation curve \textsf{G}. Together, $\textsf{SN}_0$, \textsf{G}, and $\textsf{H}_0^\pm$ indeed separate the plane into the regions III--VI.

\begin{figure*}[htp!] 
  \centering
  \includegraphics[scale=0.85]{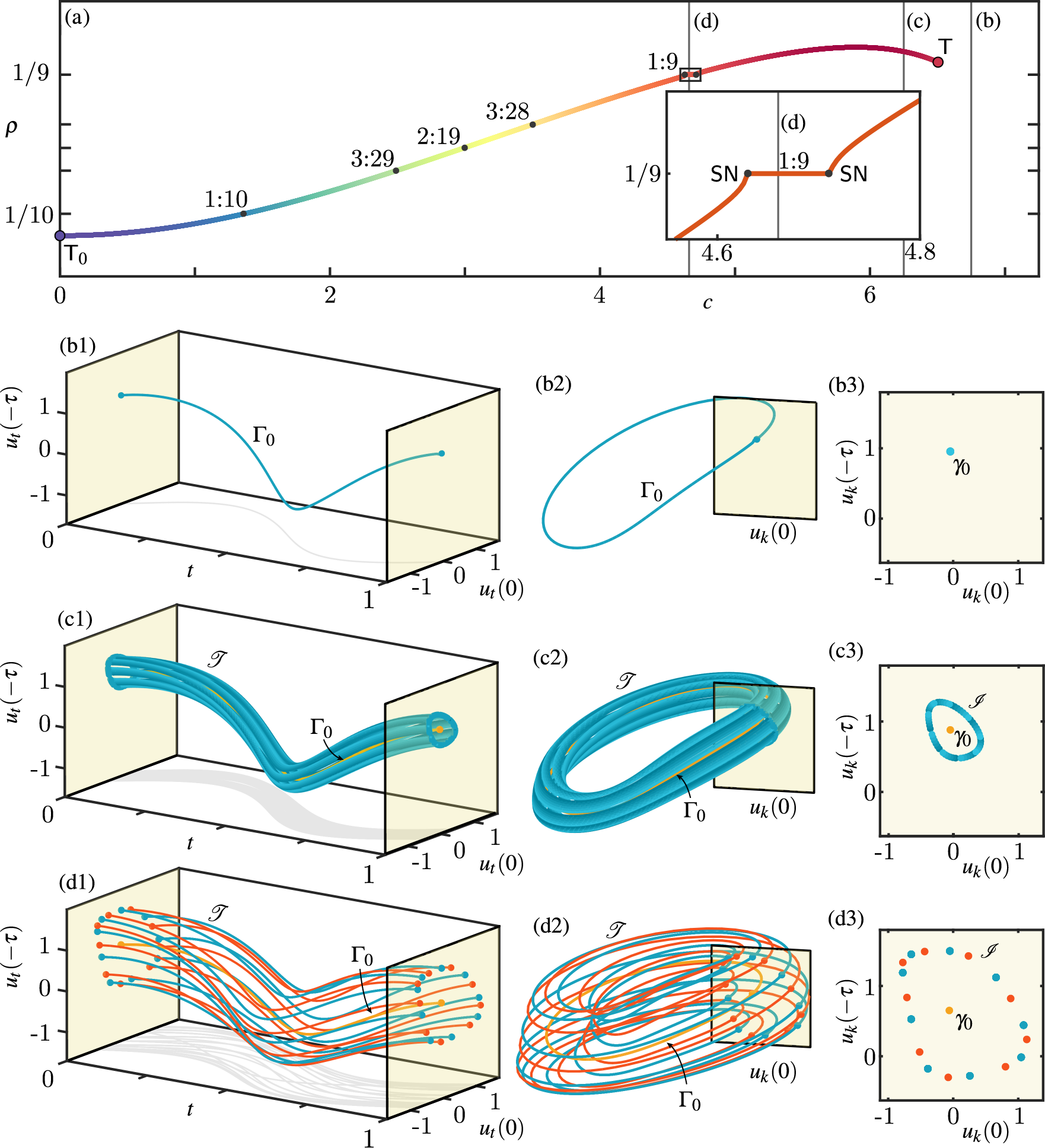}
  \caption{{One-parameter bifurcation diagram in $c$ for $\tau=2.7663$, $\alpha=0.75$ and $\omega = 2\pi$. Panel~(a) shows the graph of $\rho$, between the torus bifurcations $\textsf{T}_0$ and $\textsf{T}$, and its inset enlarges $\rho$ near the 1:9 resonance tongue. Panels~(b)–(d) illustrate representative dynamics: (b) attracting period-one orbit $\Gamma_0$ (blue) at $c=6.75$; (c) attracting quasiperiodic orbit (blue) on the invariant torus $\mathcal{T}$ with a 2-saddle period-one orbit $\Gamma_0$ (orange) at $c=6.25$; (d) attracting and 1-saddle 1:9 locked orbits (blue and red, resp.) on $\mathcal{T}$ with a 2-saddle $\Gamma_0$ (blue) at $c=4.66$. Columns~(b1)--(d1) show trajectories represented in the $(t, u_t(0), u_t(-\tau))$--space with the planes $t=0$ and $t=1$ identified, and columns~(b2)--(d2) are their toroidal embedding in $\mathbb R^3$.
Columns~(b3)--(d3) display the headpoints of the stroboscopic section (plane at $t=0$), where the period-one orbit $\Gamma_0$ is a fixed point $\gamma_0$ of the stroboscopic map $\mathcal{P}$, and the headpoints of the orbits on $\mathcal{T}$ lie on $\mathcal{I}$. }
}\label{fig:rotationumber_plot}
\end{figure*}

In regions~II and~VI, the only attracting object is the periodic orbit $\Gamma_s$. In region~II, the Hopf bifurcation curves $\textsf{H}_k$, $k\geq 1$, are subcritical: the respective saddle periodic orbit that bifurcates from the equilibrium $p_0$ exists for smaller values of $\tau$, and disappears at a saddle-node bifurcation curve $\textsf{SN}_k$. As stated in Proposition~\ref{prop:hopf}, the Hopf bifurcation curves $\textsf{H}_k$ have asymptotes $\tau=0$ as $\alpha\to\infty$ and $\alpha=1$ as $\tau\to\infty$, respectively. The associated curves $\textsf{SN}_k$ approach the same asymptote as $\alpha\to\infty$ but, at larger $\tau$, they cross the line $\alpha=1$ and enter region~VI.
In region~VI, as $\tau$ increases, one encounters the sequence $\textsf{H}_k^\pm$, $k\ge 1$, of Hopf bifurcation curves. Each as a minimum in $\tau$ and asymptotes $\alpha=1$ and $\alpha=\tfrac{1}{2}$ as $\tau\to\infty$.
Along the curve $\textsf{H}_1^\pm$ lies a
degenerate Hopf bifurcation point $\textsf{DH}_1$, at
$(\tau,\alpha)\approx(12.95,0.63)$, where the criticality of $\textsf{H}_1^\pm$ changes. From $\textsf{DH}_1$ the saddle-node
bifurcation curve $\textsf{SN}_1^\pm$ emanates; however because additional
unstable modes are present, no new attracting periodic orbit
is created. This scenario repeats along each curve
$\textsf{H}_k^\pm$, $k\geq 2$: each has a degenerate Hopf bifurcation point
$\textsf{DH}_k$ and saddle-node bifurcation curve
$\textsf{SN}_k^\pm$; these are all outside the range of Figure~\ref{fig:taualpha}.

\section{Computation of the rotation number}
\label{section:compotationofrotation}

In this section, we extend the method proposed by Bailie \emph{et al.} \cite{bailie2024resonancetongues,John}, for computing the rotation number $\rho$ of attracting invariant tori in periodically forced planar ODEs, to the case of a periodically forced DDE, specifically the pfS\&S model \eqref{eq:pss}.

\subsection{Torus dynamics in the pfS\&S model}
\label{subsection:pfSS_example}
Figure~\ref{fig:rotationumber_plot} presents a one-parameter bifurcation diagram in $c$ for fixed  $\tau = 2.7663$, $\alpha = 0.75$ and $\omega = 2\pi$. The value of $\tau$ is chosen because it corresponds to a local maximum in $c$ of the torus bifurcation curve $\textsf{T}$ in the $(\tau,c)$-parameter plane (see already Figure~\ref{fig:tau_beta_75_2pi}).
Figure~\ref{fig:rotationumber_plot}(a) shows the rotation number $\rho$ of the normally hyperbolic attracting invariant torus $\mathcal{T}$ that exists between the torus bifurcation $\textsf{T}$ at $c\approx 6.504$ and the point $\textsf{T}_0$ at $c=0$, where one finds a `trivial torus bifurcation' when $c$ is switched on. 
Selected resonance tongues are indicated by black dots on the graph of $\rho$.
The two saddle-node bifurcations of periodic orbits \textsf{SN} bounding the larger 1:9 resonance tongue were detected with \texttt{DDE-BifTool}; they are shown more clearly in the inset.
Rows~(b)--(d) show representative dynamics: a period-one orbit $\Gamma_0$, a quasiperiodic (or very long-period) orbit, and a pair of locked orbits on the invariant torus $\mathcal{T}$, respectively.
Each column displays a different projection of the dynamics.
In column~(b1)--(d1), trajectories are represented in the three-dimensional $(t, u_t(0), u_t(-\tau))$-space, with the planes $t=0$ and $t=1$.
In column~(b2)--(d2), these planes are identified to obtain an embedded torus in $\mathbb{R}^3$.
This is done by introducing $\theta = 2\pi t$ and applying the 
transformation
\begin{equation}
\label{eq:toruscoordinates}
\big((R + r\,u_t(0))\cos(\theta),\; (R + r\,u_t(0))\sin(\theta),\; r\,u_t(-\tau)\big),
\end{equation}
where $R = 1.75$ and $r = 0.75$ control the major and minor radii of the torus. 
Column~(b3)--(d3) shows the $(u_t(0),u_t(-\tau))$-plane at $t=0$ with headpoints of associated history segments. 
The period-one orbit $\Gamma_0$ defines the fixed point $\gamma_0$ of the stroboscopic map $\mathcal{P}$, and the invariant torus $\mathcal{T}$ yields the invariant circle $\mathcal{I}$.

In Figure~\ref{fig:rotationumber_plot}(a), the periodic forcing is dominating the dynamics for large $c$, and here the only attracting object is a period-one orbit $\Gamma_0$; the dynamic is seasonally dominated \cite{Tony2026}. 
When $c$ is decreased, 
the period-one orbit $\Gamma_0$ loses stability at the torus bifurcation \textsf{T}, and
an attracting invariant torus $\mathcal{T}$ is born. The torus $\mathcal{T}$ persists down to $c=0$, where it `collapses' onto the attracting periodic orbit $\Gamma_s$ of the unforced system \eqref{eq:ss} at the trivial torus bifurcation $\textsf{T}_0$.
The graph of $\rho$ is a devil's staircase—a piecewise constant and continuous, with plateaus at rational values corresponding to locked orbits (cross-sections of resonance tongues) and quasiperiodic motion in between. Each cross-section of a resonance tongue is bounded by a pair of saddle-node bifurcations of periodic orbits $\textsf{SN}$, and the largest, shown in the inset, is the 1:9 resonance tongue. The graph of $\rho$ appears smooth elsewhere only because the remaining resonance tongues are too narrow to resolve.

Figure~\ref{fig:rotationumber_plot}(b) shows the period-one orbit $\Gamma_0$ for $c=6.75$, which appears as a single string in panel~(b1), a closed loop in panel~(b2), and the fixed point $\gamma_0$ of the stroboscopic map $\mathcal{P}$ in panel~(b3).
After the torus bifurcation \textsf{T}, the attracting invariant torus $\mathcal{T}$ is born. The period-one orbit $\Gamma_0$ is now a saddle with two unstable Floquet multipliers; we say that it is a 2-saddle periodic orbit. Panel~(c) shows a quasiperiodic orbit (or of very high period) at $c=6.25$.
The orbit densely fills $\mathcal{T}$, forming a cylinder in~(c1) and a torus in~(c2), where
$\Gamma_0$ lies inside. In panel~(c3), the headpoints in the $(u_t(0),u_t(-\tau))$-plane form a smooth invariant closed curve $\mathcal{I}$ that encircles $\gamma_0$.
Row~(d) shows a pair of 1:9 locked orbits; an attracting periodic orbit and a saddle (1-saddle) periodic orbit winding around $\mathcal{T}$ at $c=4.66$. In the stroboscopic section of panel~(d3), the attracting and saddle locked orbits each correspond to headpoints lying on $\mathcal{I}$, with the period-one point $\gamma_0$ inside it.

\subsection{Computation of the rotation number}
\label{subsection:computerot}

We introduce a method to compute the rotation number $\rho$ of an attracting invariant torus $\mathcal{T}$ for periodically forced DDEs~\eqref{eqn:pDDE} over a two-parameter grid, specifically, in the $(\tau,c)$-plane for the pfS\&S model \eqref{eq:pss}.
We focus on parameter regions near torus bifurcation curves, where one expects attracting invariant tori and resonance tongues. At each grid point, we integrate \eqref{eqn:pDDE} from an initial condition and allow sufficient time for the trajectory to settle onto the attractor; we then compute $\rho$ with the procedure presented below. To ensure that we follow the same attractor as parameters vary, we use parameter sweeping in a selected parameter by increasing it step-by-step and using the last function segment of the previous computation as the  initial function for the next.

\begin{figure}[t!] 
  \centering
  \includegraphics[width=7.5cm]{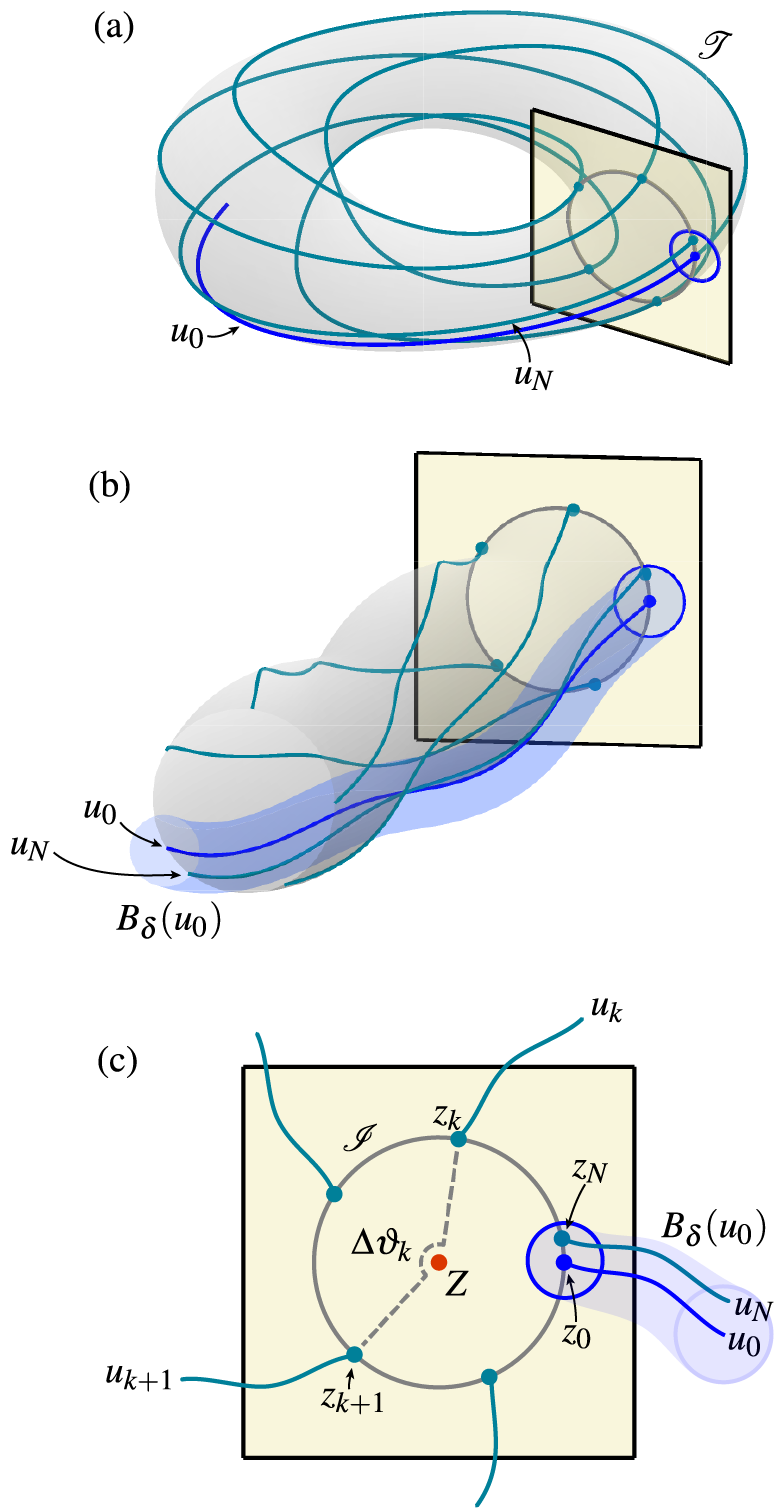}
  \caption{
Schematic of the computation of the rotation number $\rho$ for the periodically forced DDE \eqref{eqn:pDDE}. Panel (a) shows a trajectory $u(t)$ (green) winding around an attracting invariant torus $\mathcal{T}$ (grey), which intersects the stroboscopic plane (yellow). Panel~(b) illustrates how the stroboscopic map $\mathcal{P}$ creates solution segments $u_0,\dots,u_N$. The index $N$ is chosen so that $u_N \in B_\delta(u_0)$ (the blue cylinder). We show the corresponding headpoints $\{z_k\}$ in panel~(c) with their history segments, as well as the reference point ${Z}$ (red), from which we compute the oriented angular increments $\Delta\vartheta_k = \angle(z_k,{Z}, z_{k+1})$.
}\label{fig:rotnumber}
\end{figure}

To compute the rotation number $\rho$ at a given parameter point, we extend the method proposed by Bailie \emph{et al}.~ \cite{bailie2024resonancetongues,John} for periodically forced planar ODEs, to periodically forced DDEs such as the system \eqref{eq:pss}. This extension is motivated by Figure~\ref{fig:rotationumber_plot}(c3)--(d3), which shows that the headpoints of the stroboscopic sections in the $(u_t(0),u_t(-\tau))$-plane lie on a smooth invariant closed curve $\mathcal{I}$, in close analogy with the planar case.
This observation indicates that the stroboscopic map $\mathcal{P}$ acts as a circle map on $\mathcal{I}$, and this allows us to adapt the method of \cite{John}. 
The key idea is to iterate the Poincar\'{e} map until the final point is in a small neighborhood of the first point on $\mathcal{I}$.
The rotation number $\rho$ is then the average angular displacement, as in \eqref{eq:rotpractice}, over one full finite orbit of $\mathcal{P}$. 
Adapting this approach to the DDE setting faces two new
difficulties. First, the return condition cannot be imposed solely on
the headpoints in the Poincar\'{e} section, but must be formulated in the full function
space $\mathcal{C}$ to ensure that entire solution segments are close.
Second, since the dynamics evolves in an infinite-dimensional space, a
suitable projection onto $\mathbb{R}^2$ is required in order to define
and compute the angle increments.

\medskip

The method proceeds in four steps, illustrated in
Figure~\ref{fig:rotnumber}, and we describe each in turn.
We assume here that the torus $\mathcal{T}$ bifurcates from a fixed point of $\mathcal{P}$ (which can always be achieved by considering a higher iterate of $\mathcal{P}$ if necessary).

\paragraph*{Step 1: Transient integration.}
Given an initial function $\varphi\in\mathcal{C}$, and a fixed initial time $t_0=0$, we integrate \eqref{eqn:pDDE} over a transient interval $[0, T_{\mathrm{tr}}]$, where $T_{\mathrm{tr}}\in\mathbb{N}$ is chosen large enough for the trajectory $u(t)$ to settle onto the attractor; see Figure~\ref{fig:rotnumber}(a). We refer to the first history segment on $\mathcal{I}$ as $u_0$.

\paragraph*{Step 2: Stroboscopic iteration and return criterion.}
Starting from $u_0$, we integrate
\eqref{eqn:pDDE} over successive time intervals
$[T_{\mathrm{tr}}+k-1,\, T_{\mathrm{tr}}+k]$ to generate a sequence
of solution segments 
\begin{equation*}
  u_k = \mathcal{P}(u_{k-1}), \qquad k = 1, 2, \ldots
\end{equation*}
in the section. This iteration stops at the first index $N$ for which $u_N \in B_\delta(u_0)$, where $B_\delta(u_0)$ denotes
the open ball in $\mathcal{C}$ of radius $\delta > 0$ centered at
$u_0$. Measuring proximity
in the full function space $\mathcal{C}$ with the function norm $\|\cdot\|_{\mathcal{C}}$ ensures that entire solution
segments are close; see
Figure~\ref{fig:rotnumber}(b).

Three outcomes are possible.
If $N = 1$, the trajectory has converged to a fixed point of
$\mathcal{P}$, i.e., to a period-one orbit of~\eqref{eqn:pDDE}, and the
rotation number is undefined (the torus bifurcation has not yet taken place).
If $N \leq N_{\max}$, the return criterion is met: for a quasiperiodic
orbit this is guaranteed by the density of the orbit on $\mathcal{T}$ (and thus on $\mathcal{I}$); and for a $p{:}q$ locked orbit the return occurs exactly after $q$
iterates.
In either case we identify $u_N$ with $u_0$, thus, treating the orbit as
$N$-periodic and closing it accordingly.
If $N_{\max}$ is exceeded without the criterion being satisfied, we retain
$u_{N_{\max}}$ without closing the orbit and flag the instance for
post-hoc inspection; the interpretation of such cases is discussed later.

Since we do not know a priori whether the orbit is locked or 
quasiperiodic, $\delta$ and $N_{\max}$ must be chosen to accommodate both 
cases. The tolerance $\delta$ is the primary tuning parameter: it controls 
the resolution of $\rho$ by setting how tightly the orbit must close 
before $N$ is accepted. A smaller $\delta$ resolves rotation numbers with 
larger denominators $q$ and improves accuracy in the quasiperiodic case, 
but requires more iterates and, therefore, a larger $N_{\max}$. The role of 
$N_{\max}$ is purely as a safeguard: it should be large enough to ensure that the 
return criterion is met under normal operation. 
The parameters $\delta$ and $N_{\mathrm{max}}$ must be tuned jointly to the required resolution, this depends on the system and the parameter region that is investigated.

\paragraph*{Step 3: Projection onto the plane.}
The iterates $\{u_k\}$ lie on a smooth invariant curve
$\mathcal{I}_{\mathcal{C}}$ in function space.
We project them onto $\mathbb{R}^2$ via a map
$$\Pi\colon\mathcal{C}\to\mathbb{R}^2,$$ to obtain finite points
\begin{equation*}
  z_k = \Pi(u_k), \qquad k = 0, 1, \ldots, N,
\end{equation*}
which lie on the planar curve $$\mathcal{I} = \Pi(\mathcal{I}_{\mathcal{C}})\subset\mathbb{R}^2.$$
Taking the headpoints in a specified Poincar\'{e} section as in Figure~\ref{fig:rotationumber_plot}(b3)-(d3), is a natural choice for $\Pi$; see \cite{Calleja2017ResonancePI,krauskopf2003computing} for further discussion on  suitable projections for DDEs.

\paragraph*{Step 4: Angle summation.}
As reference point $Z$ in the interior of $\mathcal{I}$,
we take the componentwise mean of the projected points $\{z_k\}$ of the approximate $N$-cycle:
\begin{equation}
    \label{eq:Z}
    Z=\frac{1}{N}\sum_{k=0}^{N-1} z_k.
\end{equation}
For each $k$, we then compute the oriented angle increment
$\angle(z_k, Z, z_{k+1})$ between successive
projected points; see Figure~\ref{fig:rotnumber}(c).
The rotation number $\rho$ is then computed by the finite Birkhoff sum 
\begin{equation}
  \rho = \frac{1}{2\pi N}\sum_{k=0}^{N-1}\angle(z_k, Z, z_{k+1}).
  \label{eq:rotnumber_discrete}
\end{equation}
For periodic orbits, this finite sum is exactly the rotation number $\rho$; for quasiperiodic orbits, it approximates $\rho$ well for sufficiently large $N$.

\begin{algorithm}[H]
\label{algo}
\caption{Computation of the rotation number $\rho$}
\KwIn{
  Initial function $\varphi\in\mathcal{C}$,
  transient time $T_{\mathrm{tr}}\in\mathbb{N}$,
  tolerance $\delta>0$,
  maximum index $N_{\max}$,
  projection $\Pi\colon\mathcal{C}\to\mathbb{R}^2$}
\KwOut{Rotation number $\rho$,
Iteration reached $N$}

\BlankLine
\textbf{Step 1: Transient integration.}\\
Integrate \eqref{eqn:pDDE} with initial function $\varphi$
over $[0,T_{\mathrm{tr}}]$ to allow the trajectory to settle
onto $\mathcal{T}$. Set $u_0=u_{T_{\mathrm{tr}}}$.

\BlankLine
\textbf{Step 2: Stroboscopic iteration and return time.}\\
Set $k\leftarrow1$\;
\While{$u_k\notin B_\delta(u_0)$ \textbf{and} $k\leq N_{\max}$}{
    Compute $u_k=\mathcal{P}(u_{k-1})$\;
    $k\leftarrow k+1$\;
}
Set $N\leftarrow k$\;

\BlankLine
\uIf{$N = 1$}{
 \tcp{Period 1 orbit}
    \KwRet{}
}
\eIf{$N \leq N_{\text{max}}$}{
    Set $u_N \leftarrow u_0$ \tcp*{Period $N$ orbit}
}{
    \tcp{Orbit did not close — maximum period exceeded}
}
Store the sequence $u_0,\dots,u_N$\;

\BlankLine
\textbf{Step 3: Projection onto the plane}\\
Compute $z_k=\Pi(u_k)$ for $0\le k\le N$\;

\BlankLine
\textbf{Step 4: Angle summation.}\\
Compute $Z=\tfrac{1}{N}\sum_{k=0}^{N-1}z_k$\;
Compute $\rho=\tfrac{1}{2\pi N}\sum_{k=0}^{N-1}\angle(z_k,Z,z_{k+1})$\;
\BlankLine
\KwRet{$\rho$,\ $N$}
\end{algorithm}

The overall algorithm is presented as pseudo-code in  Algorithm~\ref{algo}.
Here, we apply it to the pfS\&S model~\eqref{eq:pss} in the $(\tau,c)$-plane with forcing frequency $\omega = 2\pi$. We fix $T_{\mathrm{tr}} = 1000$, $N_{\max} = 8000$, and $\delta = 10^{-4}$, while incrementing $\tau$ by $0.05$ and $c$ by $0.0125$. This allows us to distinguish locked $p{:}q$ orbits with denominator $q \leq N_{\max}  - 1 = 7999$. To determine the distance between history segments for the return criterion we use the supremum norm; this is advantageous compared to integral norms because it does not require evaluating a definite integral.

\medskip

\begin{figure}[t!] 
  \centering
  \includegraphics[width=10cm]{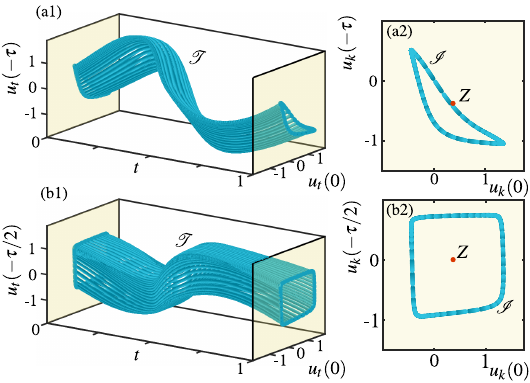}
  \caption{
  A single trajectory  on a torus $\mathcal T$ of the pfS\&S model~\eqref{eq:pss} with $\alpha = 1.25$, $c = 6.5$, $\tau = 9.11$ and $\omega = 2\pi$, shown in two different projections $\Pi_1$ and $\Pi_2$ in each row. Left column shows trajectories represented in $(t, u_t(0), u_t(-\tau))$-space (top) and $(t, u_t(0), u_t(-\tau/2))$-space  (bottom), respectively.
Right column shows the corresponding projected stroboscopic sections at $t=0$ and the reference point $Z$ with the invariant curve $\mathcal I$. 
}\label{fig:rotnumberissue}
\end{figure}

\begin{figure*}[t!] 
  \centering
    \includegraphics[scale=.87]{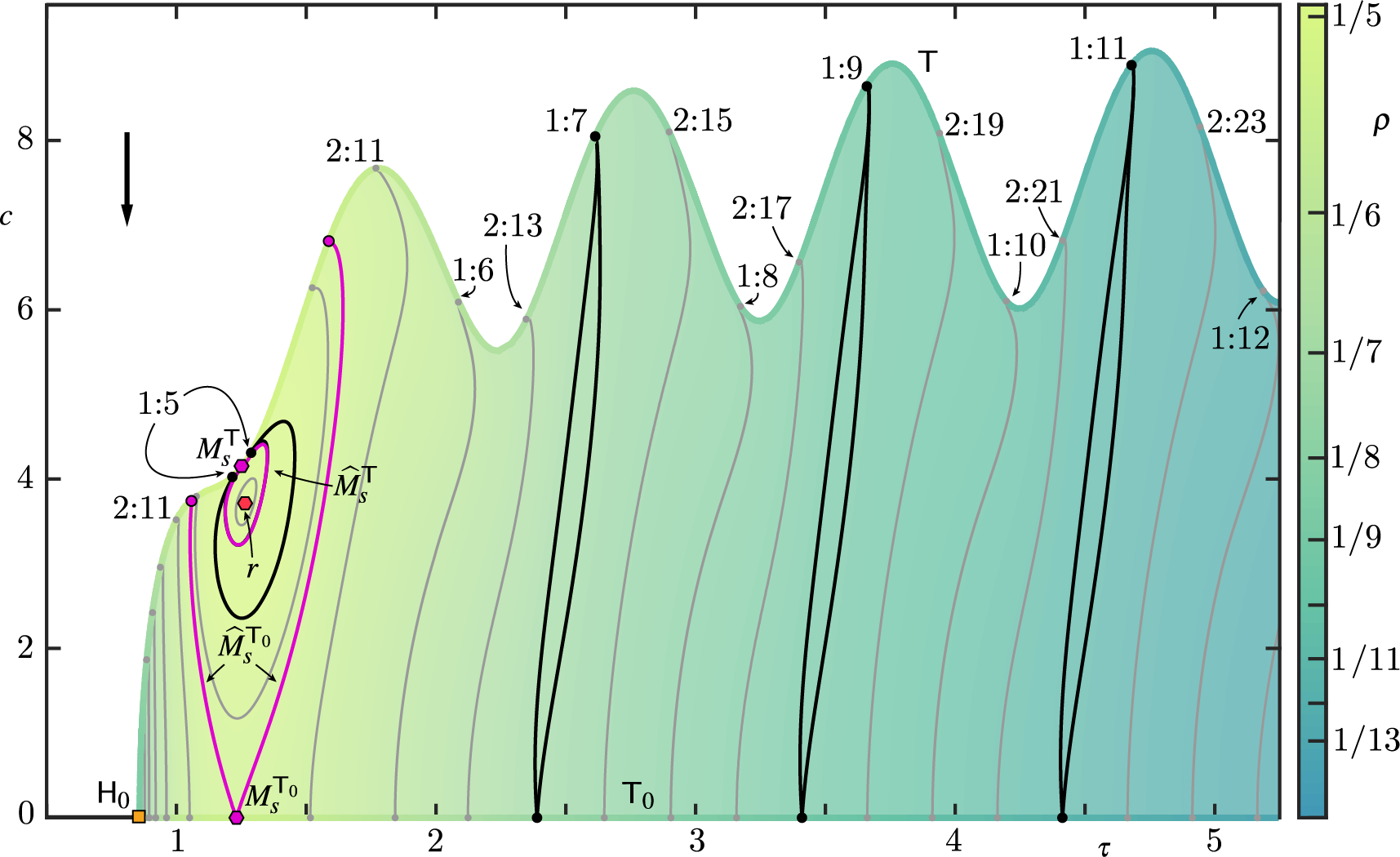}
\caption{Bifurcation set and rotation number $\rho$ of the pfS\&S model~\eqref{eq:pss} in the $(\tau,c)$-plane for $\alpha = 1.25$ and $\omega = 2\pi$. In the white region, the period-one orbit $\Gamma_0$ is attracting, while in the coloured region the invariant torus $\mathcal{T}$ is attracting, with the rotation number $\rho=\rho(\tau,c)$ indicated by the colorbar. This region is bounded by the torus bifurcation curves $\textsf{T}$ and $\textsf{T}_0$, which emanate from the Hopf bifurcation $\textsf{H}_0$. We show $p{:}q$ resonance tongues and their root points of the form $1{:}q$ and $2{:}q$ in black and grey, as well as critical points of $\rho$. Boundary max-saddle points ${M}_s^{\textsf{T}}$ and ${M}_s^{\textsf{T}_0}$, and their associated separatrices $\widehat{M}_s^{\textsf{T}}$ and $\widehat{M}_s^{\textsf{T}_0}$, are shown in pink, and the interior maximum $r$ in red.
}\label{fig:tau_beta_125_2pi}
\end{figure*}

The validity of Steps~3 and 4 depends critically on the choice of 
projection $\Pi\colon\mathcal{C}\to\mathbb{R}^2$, which must map the
invariant curve $\mathcal{I}_{\mathcal{C}}$ to a simple closed curve
$\mathcal{I}$ in $\mathbb{R}^2$, with the reference point $Z$ from \eqref{eq:Z} in its interior. 
In general, $\Pi(\mathcal{I}_\mathcal{C})$ may not be an embedding of $\mathbb{S}^1$ in $\mathbb{R}^2$.
On the other hand, Takens' embedding 
theorem~ \cite{takens2006detecting,sauer1991embedology} guarantees that a generic map using three delay-evaluation coordinates results in the original invariant circle $\mathcal{I}_{\mathcal{C}} \subset \mathcal{C}$ being represented by an embedding of $\mathbb S^1$ into $\mathbb{R}^3$; specifically, this is (generically) the case for
\[
\varphi \mapsto \bigl(\varphi(0),\;\varphi(-\tau/2),\;\varphi(-\tau)\bigr).
\]
Because Step~4 of Algorithm~\ref{algo} 
requires a two-dimensional setting to compute angles, one must then 
project from $\mathbb{R}^3$ onto $\mathbb{R}^2$. As a practical 
choice, we consider the two projections 
\begin{equation}
    \label{eq:Pi1}
    \Pi_1(\varphi)=(\varphi(0),\varphi(-\tau))
\end{equation} (the usual headpoints) and 
\begin{equation}
    \label{eq:Pi2}
    \Pi_2(\varphi)=(\varphi(0),\varphi(-\tau/2))
\end{equation}
(half-delay headpoints). The expectation is that at least one yields a
well-shaped $\mathcal{I}$ in $\mathbb{R}^2$ with $Z$ in its interior; as we will see, this approach works well for the pfS\&S model \eqref{eq:pss}.

Past a torus bifurcation, the invariant circle $\mathcal{I}_\mathcal{C}$
is small and the projection $\Pi$ onto the headpoints of the stroboscopic section produces a simple and convex curve $\mathcal{I}$.
This is
illustrated in Figure~\ref{fig:rotationumber_plot}(c3) and (d3) for the
standard projection $\Pi_1$ in the $(u_t(0),u_t(-\tau))$-plane. As parameters move
away from the bifurcation, the torus grows and may distort, and the
projected curve $\mathcal{I}$ can become problematic; for instance,
\cite{Calleja2017ResonancePI} reports a case where
a projected invariant curve self-intersects to form a figure-eight. 
Figure~\ref{fig:rotnumberissue} provides an example where the head-point projection $\Pi_1$ in row~(a) fails, while $\Pi_2$ of the half-delay headpoints in row~(b) yields a suitable projection. The left column~(a1)--(b1) shows the trajectory in the
$(t,u_t(0),u_t(-\tau))$-space and $(t,u_t(0),u_t(-\tau/2))$-space,
respectively, from $t=0$ to $t=1$; and the right column~(a2)--(b2) shows the corresponding
headpoints forming $\mathcal{I}$ with the centroid $Z$. 
While both curves are simple, the point $Z$ computed with \eqref{eq:Z} lies outside the banana-shaped $\Pi_1(\mathcal{I}_\mathcal{C})$ in panel~(a2). As a result, the infinite Birkhoff sum is zero, and we get a spurious value of $\rho = 6.59\times10^{-5}$. However, $\Pi_2(\mathcal{I}_\mathcal{C})$ in panel~(b2) produces a simple closed curve with $Z$ well inside it, and we obtain the value $\rho = 0.058$.

Since it is not possible to determine beforehand where the standard projection $\Pi_1$ fails, we produce $\rho$ for both projections $\Pi_1$ and $\Pi_2$. At the end of the computation, we perform a diagnostic check of the continuity of $\rho$. If the continuity diagnostic reveals a jump discontinuity in $\rho$ under one projection but not the other, this indicates a problem with the projection: the reference point $Z$ may have drifted outside $\mathcal{I}$. We then discard the projection that returns an approximate value of 0 and keep the other one. This ensures a consistent computation of $\rho$ throughout the parameter region, as shown in Figures~\ref{fig:tau_beta_125_2pi} and \ref{fig:tau_beta_75_2pi}.

The computation of the rotation number $\rho$ of a DDE in parameter space is expensive. For instance, the grid in the $(\tau,c)$-plane for the pfS\&S model \eqref{eq:pss} shown in Figure~\ref{fig:tau_beta_125_2pi} requires $768{,}000$ rotation number computations in total, taking approximately $8$ hours with the third-order Adams--Bashforth method in \texttt{Julia} without parallelization on an Apple M1 (16\,GB).

\section{Bifurcation analysis of the pfS\&S model}
\label{section:bif_pss}
We now analyze the torus dynamics of the pfS\&S model~\eqref{eq:pss} in the $(\tau,c)$-plane for forcing frequency $\omega = 2\pi$. To this end, we first consider in Section~\ref{section:bifalpha125} the simpler case that $\alpha = 1.25$. We then present in Section~\ref{section:bifanal075} the resonance structure for the richer and more physically relevant case $\alpha = 0.75$.

\subsection{Case $\alpha = 1.25$}
\label{section:bifalpha125}
Figure~\ref{fig:tau_beta_125_2pi} shows the bifurcation diagram of Eq.~\eqref{eq:pss} in the $(\tau,c)$-plane for $\omega = 2\pi$ and $\alpha = 1.25$. 
 The torus bifurcation curve $\textsf{T}$ emerges 
from the supcritical Hopf bifurcation $\textsf{H}_0$ near $(\tau,c)\approx(0.8580,0)$, which is the Hopf bifurcation $\textsf{H}_0$ of the unforced system~\eqref{eq:ss} from Proposition~\ref{prop:hopf}(iv). The curve \textsf{T} separates the $(\tau,c)$-plane into two regions. Above $\textsf{T}$, the seasonal forcing dominates and the period-one orbit $\Gamma_0$ is attracting; below it, $\Gamma_0$ is a saddle period-one orbit and an attracting invariant torus $\mathcal{T}$ exists. This region with an attracting torus is bounded by the torus bifurcation curves $\textsf{T}$ and $\textsf{T}_0$ (with $c=0$), and it extends seemingly periodically to arbitrarily large values of $\tau$.
The rotation number $\rho$ in this region is computed with Algorithm~\ref{algo} with downward sweeps in $c$ for fixed values of $\tau$, and it is illustrated in color according to the colorbar. 
Along the torus bifurcation curves $\textsf{T}$ and $\textsf{T}_0$, the torus does not yet exist, but the rotation number $\rho$ is still well-defined. At $\textsf{T}$, the critical Floquet multipliers are 
$e^{\pm 2\pi i \kappa}$, and $\rho=\kappa$ (or $\rho = 1-\kappa$) serves naturally as the rotation number. At $\textsf{T}_0$, the rotation number is defined as the ratio of the forcing period and the natural period $T$ of the autonomous periodic orbit $\Gamma_s$ from Eq.~\eqref{eq:ss}; since the forcing frequency is $\omega = 2\pi$ (period $1$), this gives $\rho = 1/T$. Importantly, $\rho$ varies smoothly along the torus bifurcation curves $\textsf{T}$ and $\textsf{T}_0$, allowing us to locate resonance points at its rational values. In Figure~\ref{fig:tau_beta_125_2pi}, we marked with dots all $p{:}q$ resonance points of the form $1{:}q$ and $2{:}q$, of which there are finitely many because the rotation number $\rho$ is bounded away from zero; these serve as root points of the associated shown resonance tongues. The boundaries of 
larger resonance tongues (black) are obtained with \texttt{DDE-BifTool} by tracking saddle-node 
bifurcations $\textsf{SN}$ of the pair of locked orbits. The remaining resonance tongues (grey) are 
found as level sets of $\rho$. 

\begin{figure*}[t!] 
  \centering
    \includegraphics[scale=0.86]{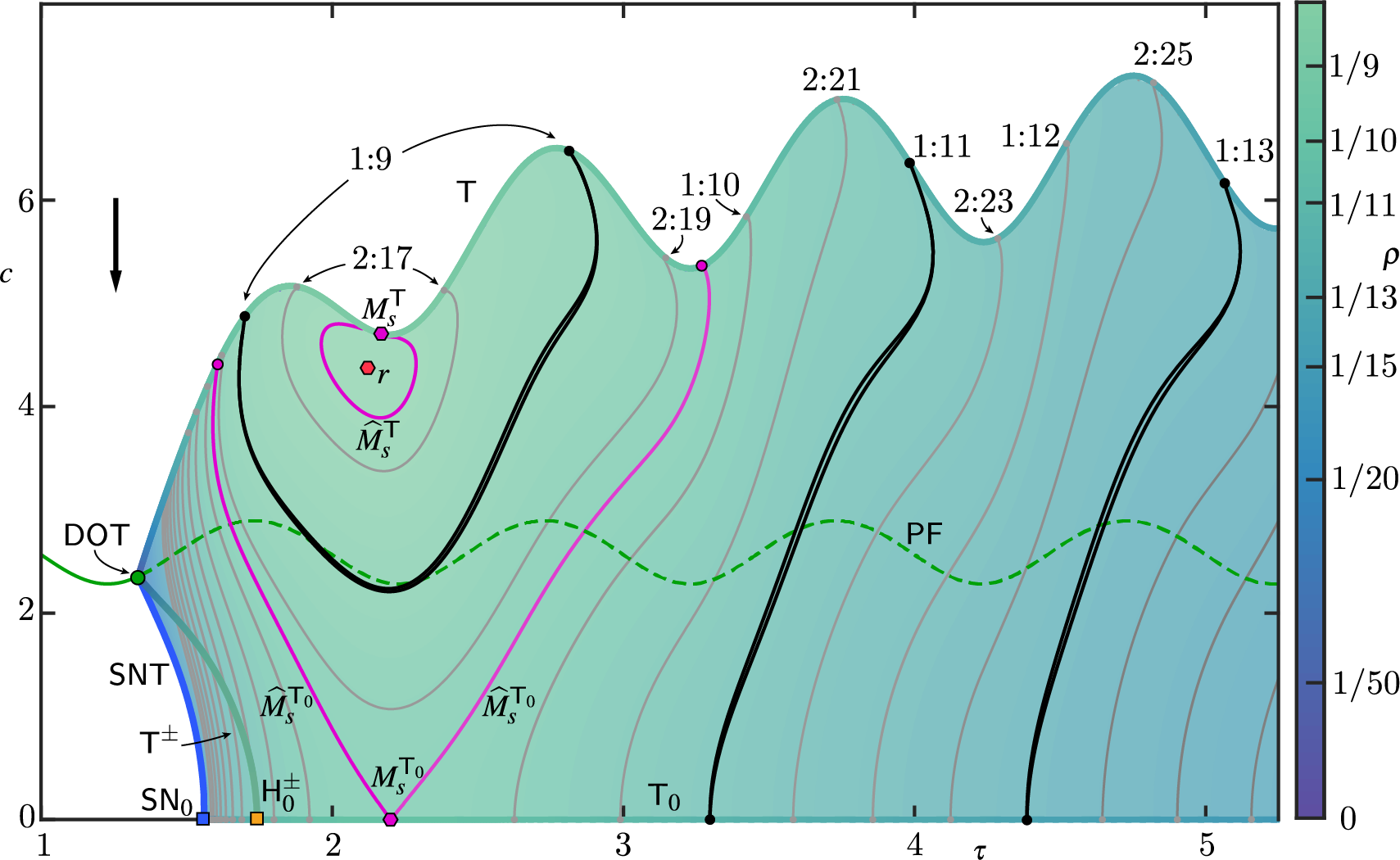}
\caption{
Bifurcation set and rotation number $\rho$ of the pfS\&S model~\eqref{eq:pss} in the $(\tau,c)$-plane for $\alpha = 0.75$ and $\omega = 2\pi$. In the white region, above the pitchfork of period-one orbits $\textsf{PF}$ in green, the period-one orbit $\Gamma_0$ is attracting, and the period-one orbits $\Gamma^{\pm}$ are attracting below it. In the coloured region, the invariant torus $\mathcal{T}$ is attracting, with the rotation number $\rho$ indicated by the colorbar, which reaches zero at the double-one bifurcation of tori point $\textsf{DOT}$. This region is bounded by the torus bifurcation curves $\textsf{T}$ and $\textsf{T}_0$ and the locus of $\textsf{SNT}$. We show $p{:}q$ resonance tongues and their root points of the form $1{:}q$ and $2{:}q$ with $q\leq 50$ in black and grey, as well as critical points of $\rho$: boundary max-saddle points ${M}_s^{\textsf{T}}$ and ${M}_s^{\textsf{T}_0}$, and their associated separatrices $\widehat{M}_s^{\textsf{T}}$ and $\widehat{M}_s^{\textsf{T}_0}$, are shown in pink, and the interior maximum $r$ in red. Along $c=0$, we show the saddle-node of periodic orbit $\textsf{SN}_0$, the gluing bifurcation $\textsf{G}$, and the subcritical Hopf bifurcations $\textsf{H}_0^{\pm}$.
}\label{fig:tau_beta_75_2pi}
\end{figure*}

In Figure~\ref{fig:tau_beta_125_2pi}, the resonance tongues do not overlap; that is, an attracting invariant and normally hyperbolic torus $\mathcal{T}$ exists throughout the colored region. Hence, the graph of $\rho=\rho(\tau,c)$ is defined and continuous everywhere in the region, and it has the following properties: for every rational value of $\rho$, there exists a plateau of locked orbits -- a resonance tongue -- containing an attracting and a saddle periodic orbit. Quasiperiodic motion occurs at irrational values of $\rho$ between the plateaus. Thus, the graph of $\rho$ is piecewise constant, reminiscent of a devil's staircase but extended to two parameters, for which the name \emph{devil's terraces} has been coined~ \cite{JIN1996442,John}. 
For all rational values of the rotation number $\rho$, we find root 
points of resonance tongues on the torus bifurcation curve 
$\textsf{T}_0$, all of which connect to corresponding root points on the torus bifurcation curve 
$\textsf{T}$. On $\textsf{T}_0$, the supercritical Hopf point 
$\textsf{H}_0$ is a local minimum of $\rho$, from which $\rho$ 
increases as $\tau$ increases until reaching a local maximum 
$M_s^{\textsf{T}_0}$; for larger $\tau$, the rotation number $\rho$ decreases again. This local maximum $M_s^{\textsf{T}_0}$
plays a special role: it is associated with two separatrix level sets 
$\widehat{M}_s^{\textsf{T}_0}$ along which $\rho$ is constant, each 
connecting to $\textsf{T}$ at the corresponding value of $\rho$. Along $\textsf{T}$, we similarly 
identify a local maximum $M_s^{\textsf{T}}$, where the value of $\rho$ 
at $M_s^{\textsf{T}}$ is larger than at $M_s^{\textsf{T}_0}$. 
The separatrix $\widehat{M}_s^{\textsf{T}}$connects to itself at the point ${M}_s^{\textsf{T}}$ to form a closed loop.

The organization of resonance tongues around extrema of $\rho$ on torus 
bifurcation curves was studied by Terrien \emph{et al.}~ \cite{terrien}. Building on this, Bailie \emph{et al.}~ \cite{bailie2024resonancetongues,John}, working with a periodically forced planar ODE, extended this picture into the interior of the region of attracting invariant tori by applying a Morse-theory inspired framework for determining the structure of the graph of $\rho$. In this framework, critical points of $\rho$ — both on the torus bifurcation curves \textsf{T} and $\textsf{T}_0$ and in the interior of the resonance region — are classified as maxima, minima, or saddles. The local maxima $M_s^{\textsf{T}}$ and $M_s^{\textsf{T}_0}$ of $\rho$ on the torus bifurcation curves act as boundary saddles of the resonance 
region~ \cite{John}, and their separatrices $\widehat{M}_s^{\textsf{T}}$ and $\widehat{M}_s^{\textsf{T}_0}$ partition the resonance 
region, grouping resonance tongues into families according to which side of the separatrix they lie on. In Figure~\ref{fig:tau_beta_125_2pi}, we identify an interior maximum $r$ of $\rho$, located away 
from the torus bifurcation curves: it is a plateau if $\rho(r)$ is 
rational, or an isolated point if $\rho(r)$ is irrational. Overall, the critical points of $\rho$ organizes the resonance tongues into four families. Resonance tongues contained within the region enclosed by $\widehat{M}_s^{\textsf{T}}$ form annulus-like structures surrounding the interior maximum $r$; in particular,  they are not rooted to any torus bifurcation 
curve. Resonance tongues outside this region but within the domain bounded by $\widehat{M}_s^{\textsf{T}_0}$ and $\textsf{T}$, such as the $1$:$5$ resonance 
tongue, possess two root points on $\textsf{T}$, one on each side of $M_s^{\textsf{T}}$. Finally, all resonance tongues lying outside this domain entirely, such as the $2$:$11$ resonance tongue, are attached to $\textsf{T}$ on one side and to $\textsf{T}_0$ on the other; these are distinguished as two families, one on either side of $\widehat{M}_s^{\textsf{T}_0}$, as discussed above.

For sufficiently large $\tau$ (above 2), the torus bifurcation curve $\textsf{T}$ oscillates in $c$ with period one in $\tau$, and seems to converge to a sine wave with $c$ ranging between about $6$ and $9$. Moreover, the entire resonance structure also appears to repeat with period one in $\tau$: for example, the $1{:}7$ resonance tongue closely resembles the $1{:}9$, $1{:}11$, and higher odd-denominator tongues, each shifted by one in $\tau$. This self-similar reappearance mechanism is studied in detail in Section~\ref{section:global}.

\subsection{Case $\alpha = 0.75$}
\label{section:bifanal075}

We now turn to the more physically relevant and more complex case $\alpha = 0.75$.
Figure~\ref{fig:tau_beta_75_2pi} shows the two-parameter bifurcation 
diagram of Eq.~\eqref{eq:pss} in the $(\tau,c)$-plane for fixed 
$\omega = 2\pi$ and $\alpha = 0.75$ with the graph of $\rho$ represented again by color according to the colorbar. The first impression is that the resonance structure is very similar to that of Figure~\ref{fig:tau_beta_125_2pi}. 
There are torus bifurcation curves $\textsf{T}$ and $\textsf{T}_0$ with boundary saddles ${M}_s^{\textsf{T}}$ and ${M}_s^{\textsf{T}_0}$, their associated separatrices $\widehat{M}^{\textsf{T}}$ and $\widehat{M}^{\textsf{T}_0}$, and the interior maximum $r$; together, these govern the overall resonance structure shown in Figure~\ref{fig:tau_beta_125_2pi}. 
For large delays $\tau$, the resonance structure also repeats, as we will discuss in Section~\ref{section:global}.

However, the overall dynamical landscape in Figure~\ref{fig:tau_beta_75_2pi} is considerably richer, because there are additional bifurcation curves and codimension-two points for $\alpha = 0.75$: the 
pitchfork bifurcation of period-one orbits $\textsf{PF}$, a pair of 
torus bifurcation curves $\textsf{T}^\pm$, 
and the double-one torus  bifurcation point $\textsf{DOT}$. Along the zero-forcing line ($c=0$), we also indicate the subcritical Hopf bifurcations $\textsf{H}_0^\pm$, the gluing bifurcation $\textsf{G}$, and the fold of periodic orbits $\textsf{SN}_0$ of the autonomous system \eqref{eq:ss}. Note that the value of the rotation number $\rho$ at the double-one torus bifurcation point $\textsf{DOT}$ reaches zero. This is why we show the $p$:$q$ resonance tongues of the form $1$:$q$ and $2$:$q$ again, but now only for $q\leq 50$; notice how the shown resonance tongues accumulate near the point \textsf{DOT}.

In Figure~\ref{fig:tau_beta_75_2pi}, the pitchfork bifurcation 
$\textsf{PF}$ separates the $(\tau,c)$-plane into two regions. Above 
$\textsf{PF}$, there exists a single symmetric attracting periodic  
orbit $\Gamma_0$ of period $T=1$.
Below $\textsf{PF}$, there exists a 
pair of symmetry-related periodic orbits $\Gamma^\pm$ (also of period $T=1$) and $\Gamma_0$ is now unstable. At $(\tau,c)\approx(1.344,2.370)$ on $\textsf{PF}$ lies the double-one 
bifurcation of tori point $\textsf{DOT}$, where two nontrivial Floquet 
multipliers are equal to $+1$. This is a codimension-two bifurcation 
whose unfolding at 
the level of a Poincar\'e section is that of the $1$:$2$ strong resonance~ \cite{Kuznetsov}. From $\textsf{DOT}$ 
emanate three torus bifurcation curves: the curve $\textsf{T}$ above 
$\textsf{PF}$, where $\Gamma_0$ loses stability and an attracting 
invariant torus $\mathcal{T}$ is born, and the curve 
$\textsf{T}^\pm$ below $\textsf{PF}$, where $\Gamma^\pm$ lose stability 
to create a pair of saddle-type invariant tori.
The torus bifurcation curve $\textsf{T}$ behaves analogously to the case 
$\alpha = 1.25$: as $\tau$ increases, it converges to a nearly sinusoidal 
shape with period one in $\tau$. The curve $\textsf{T}^\pm$, on the 
other hand, connects to the zero-forcing line, at the subcritical Hopf bifurcation 
point $\textsf{H}_0^\pm$. After a 
`gluing bifurcation' of the saddle-type invariant tori, the resulting single saddle torus merges with the attracting torus $\mathcal{T}$ at a `saddle-node 
of tori' $\textsf{SNT}$; see Section~\ref{subsubsection:bistability} 
for details. Crucially for now, the region with  attracting invariant torus  $\mathcal{T}$ in Figure~\ref{fig:tau_beta_75_2pi} is bounded by the torus bifurcation curves $\textsf{T}$ and $\textsf{T}_0$, and the locus $\textsf{SNT}$.

\begin{figure}[t!] 
  \centering
  \includegraphics[width=9cm]{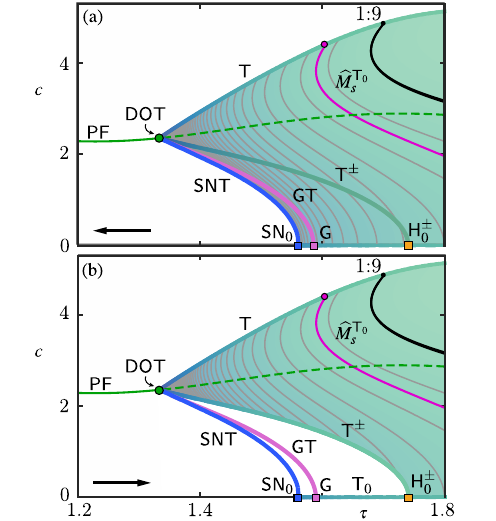}
  \caption{
  Enlarged view of the bifurcation diagram from Figure~\ref{fig:tau_beta_75_2pi} showing the region of bistability bounded between the locus \textsf{SNT} of the saddle-node of tori and torus bifurcation curves $\textsf{T}^\pm$ and $\textsf{T}_0$. The rotation number $\rho$ is given according to the color scheme shown in Figure~\ref{fig:tau_beta_75_2pi} as $c$ is fixed and $\tau$ is swept to the left in (a) and to the right in (b). Also shown is the locus \textsf{GT} of the gluing bifurcation of tori.}\label{fig:tau_c_Bistability}
\end{figure}

\begin{figure}[t!] 
  \centering
  \includegraphics[width=7.35cm]{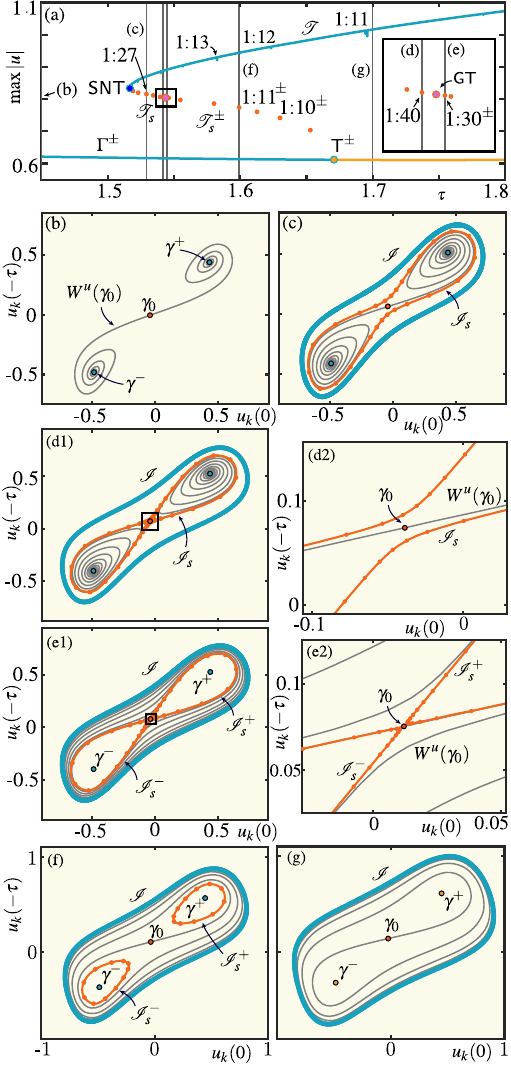}
  \caption{
One-parameter bifurcation diagram in $\tau$ for $\alpha = 0.75$, $\omega = 2\pi$ and $c=1$. Attracting objects are in blue, 1-saddles in yellow, and 2-saddles in orange. Panel (a) shows the bifurcation diagram in $\tau$ in terms of the max-norm. The attracting period-one orbits $\Gamma^\pm$ and attracting invariant torus $\mathcal T$ are shown, and the branches of saddle-type invariant tori $\mathcal T_s$ and $\mathcal T_s^\pm$ are represented by selected locked orbits. The inset zooms in on the gluing bifurcation $\textsf{GT}$.
Remaining panels show the projected headpoints of the relevant invariant objects in the stroboscopic section for   $\tau=1.45$ (b), $1.5292$ (c),  $1.5418$ (d1), $1.5451$ (e1), $1.5992$ (f), $1.7$ (g), where a period-one orbit $\Gamma$ becomes a fixed point $\gamma$, and the headpoints of a locked-orbit on a torus $\mathcal{T}$ lie on invariant circle $\mathcal{I}$. Panels~(d2)-(e2) are enlargement of panels~(d1)-(e1). Also shown is the unstable manifold $W^u(\gamma_0)$ (grey). 
}\label{fig:bistability}
\end{figure}

\subsubsection{Bistability of tori for small $\tau$}\label{subsubsection:bistability}
Figure~\ref{fig:tau_c_Bistability} shows an enlargement of Figure~\ref{fig:tau_beta_75_2pi} in the region $(\tau,c)\in[1.2,1.8]\times[0,5]$, centered around the double-one bifurcation point  $\textsf{DOT}$. In panels~(a) and (b), the values of the rotation number $\rho$ are obtained by sweeping $\tau$ in opposite directions for fixed values of $c$; this reveals a region of bistability between the attracting invariant torus $\mathcal{T}$ 
and the period-one orbits $\Gamma^\pm$. In addition to the torus bifurcation curves $\textsf{T}$ and $\textsf{T}^\pm$, it shows bifurcation loci  emanating from the bifurcation point $\textsf{DOT}$: the locus $\textsf{SNT}$ of `saddle-nodes of tori', which connects to the point $\textsf{SN}_0$ at $c=0$, and the locus \textsf{GT} of `gluing bifurcation of tori', which connects to the point $\textsf{G}$ at $c=0$.

The bistability region in Figure~\ref{fig:bistability} is 
bounded by the locus $\textsf{SNT}$, and the curves $\textsf{T}_0$ and $\textsf{T}^\pm$, and 
contains the locus $\textsf{GT}$ in its interior. Near $\textsf{SNT}$, the attracting invariant torus $\mathcal{T}$ and a saddle-type invariant torus $\mathcal{T}_s$ `collide' and annihilate. Near 
$\textsf{GT}$, a pair of saddle-type tori $\mathcal{T}_s^\pm$ `merge' into a single saddle-type torus $\mathcal{T}_s$. We stress that these loci are 
not smooth bifurcation curves: theory predicts that 
the collision or merging of tori is preceded by a loss of normal 
hyperbolicity and a sequence of torus breakdown events involving global 
bifurcations such as homoclinic tangencies; see, e.g., \cite{KKChenciner}.

\smallskip
We now look at the transition through the bistability region more closely. 
Figure~\ref{fig:bistability} shows a one-parameter bifurcation diagram in $\tau$ for $c=1$, that is, along the corresponding slice of Figure~\ref{fig:tau_c_Bistability}.
Figure~\ref{fig:bistability}(a) presents, in terms of the maximum of the solution $u(t)$, the period-one orbits $\Gamma^\pm$ and the invariant torus $\mathcal T$. 
The stable branch of $\mathcal T$ above the point \textsf{SNT} was found by integration with parameter sweeping; note that large resonance tongues appear as `local dips' because the locked solution does not attain the maximum value of $u(t)$ on the torus $\mathcal T$; see \cite{KKGZT}. The branch of saddle tori between \textsf{SNT} and $\textsf{T}^\pm$ cannot be found by forward integration; it is represented here by selected $p{:}q$ locked saddle periodic orbits. Each of these periodic orbits was continued in the $(\tau,c)$-plane from the zero-forcing line up to $c=1$. The inset shows an enlargement near  \textsf{GT}, where 
the symmetric pair of invariant tori $\mathcal{T}_s^\pm$ `merges' to form a single torus $\mathcal{T}_s$. 
The remaining panels show representative phase portraits of the stroboscopic Poincar\'e map as the headpoints in the $(u_t(0), u_t(-\tau))$-plane. Here, the period-one orbits $\Gamma_0$ and $\Gamma^\pm$ are represented by the fixed points $\gamma_0$ and $\gamma^\pm$, and the invariant torus $\mathcal{T}$, and the saddle-type tori $\mathcal{T}_s$ and $\mathcal{T}^\pm_s$, by the invariant circles $\mathcal I$, $\mathcal I_s$ and $\mathcal I_s^\pm$ respectively. The invariant circles $\mathcal I_s^\pm$ and $\mathcal I _s$ are found by computing the respective locked orbit on them and then completing the invariant curve via spline interpolation.
We also plot the one-dimensional unstable manifold $W^u(\gamma_0)$. 

For small values of $\tau$, as in Figure~\ref{fig:bistability}(b), the fixed point $\gamma_0$ is a saddle and $\gamma^\pm$ are attracting. The unstable manifold $W^u(\gamma_0)$ has two branches that converge to $\gamma^+$ and $\gamma^-$, respectively. 
As $\tau$ increases past the saddle-node of tori \textsf{SNT} at $\tau\approx 1.5175$, the system now has an attracting torus $\mathcal{T}$ and a saddle-type torus $\mathcal{T}_s$. The resulting phase portraits of $\mathcal P$ is shown in panel~(c), with the outer attracting invariant circle $\mathcal I$ and the inner saddle invariant circle $\mathcal I_s$. Note that the dynamics on $\mathcal T$ is quasiperiodic (or of very high period), while on $\mathcal T_s$ there is a 1:27 locked orbit.
At $\tau = 1.5418$, we obtain a topologically equivalent scenario 
to panel~(c), but with a $1$:$40$ locked orbit on $\mathcal{T}_s$. The 
associated saddle invariant circle $\mathcal{I}_s$, shown in panel~(d1), approaches 
the fixed point $\gamma_0$ closely; this is seen more clearly in the 
enlargement panel~(d2). As $\tau$ is increased past $\textsf{GT}$, the single saddle-type torus 
$\mathcal{T}_s$ splits into a conjugate pair $\mathcal{T}_s^-$ and 
$\mathcal{T}_s^+$. The phase portraits of $\mathcal P$ at 
$\tau = 1.5451$, shown in panel~(e1), displays the associated invariant 
circles $\mathcal{I}_s^+$ and $\mathcal{I}_s^-$, each represented by a 
$1$:$30$ locked orbit; note that $\mathcal I_s^\pm$ come extremely close to the fixed point 
$\gamma_0$, even in the enlargement panel~(e2).
As $\tau$ is increased further to $\tau =1.5992$, as in panel~(f), the saddle 
invariant circles $\mathcal{I}_s^\pm$, now represented by $1$:$11$ locked 
orbits, have shrunk sufficiently to be clearly distinguished from one 
another.
For increasing $\tau$, the saddle-type tori $\mathcal{T}_s^\pm$ continue 
to shrink, and they eventually disappear at the pair of subcritical torus bifurcations 
$\textsf{T}^\pm$ at $\tau \approx 1.672$. Beyond this point, 
illustrated in panel~(g) for $\tau = 1.7$, the periodic orbits $\Gamma^\pm$ are 
of saddle type and the attracting torus $\mathcal{T}$, whose invariant 
circle $\mathcal{I}$ remains densely filled, is globally attracting.

\smallskip

We conclude that the one-parameter bifurcation diagram of 
Figure~\ref{fig:bistability} closely mirrors the transitions found in 
the unforced S\&S model~\eqref{eq:ss} shown in 
Figure~\ref{fig:oneparameteralpha075}. 
Observe that the double-one torus bifurcation $\textsf{DOT}$ acts as an organizing 
center for the bifurcations of tori in Figure~\ref{fig:tau_c_Bistability}, 
in direct analogy with $\textsf{DZ}$ organizing the bifurcations of 
periodic orbits of~\eqref{eq:ss} in Figure~\ref{fig:taualpha}. Under 
periodic forcing, equilibria $p_0$ and $p^\pm$ lift to period-one orbits 
$\Gamma_0$ and $\Gamma^\pm$, and their bifurcations lift accordingly: 
bifurcations of equilibria become bifurcations of periodic orbits, and 
bifurcations of periodic orbits become bifurcations of invariant tori. 
This correspondence is summarized in 
Table~\ref{tab:bifurcation_analogy}.
\begin{table}[!htb]
\caption{\label{tab:bifurcation_analogy}Analogy between bifurcations of the autonomous system \eqref{eq:ss} and the periodically forced system \eqref{eq:pss}.}
\begin{tabular}{cc}
\hline
\textbf{autonomous} & \textbf{periodically forced} \\
\hline
Pitchfork of equilibria $\textsf{P}$ & Pitchfork of period-1 orbits $\textsf{PF}$ \\
Hopf bifurcations $\textsf{H}_0, \textsf{H}_0^\pm$ & Torus bifurcations $\textsf{T}, \textsf{T}^\pm$ \\
Fold of periodic orbits $\textsf{SN}_0$ & `Fold of tori' $\textsf{SNT}$ \\
gluing of periodic orbits $\textsf{G}$ & `gluing of tori' $\textsf{GT}$ \\
Double-zero bifurcation $\textsf{DZ}$ & Double-one bifurcation of tori $\textsf{DOT}$ \\
\hline
\end{tabular}
\end{table}

\begin{figure}[t!] 
  \centering
    \includegraphics[width=12.85cm]{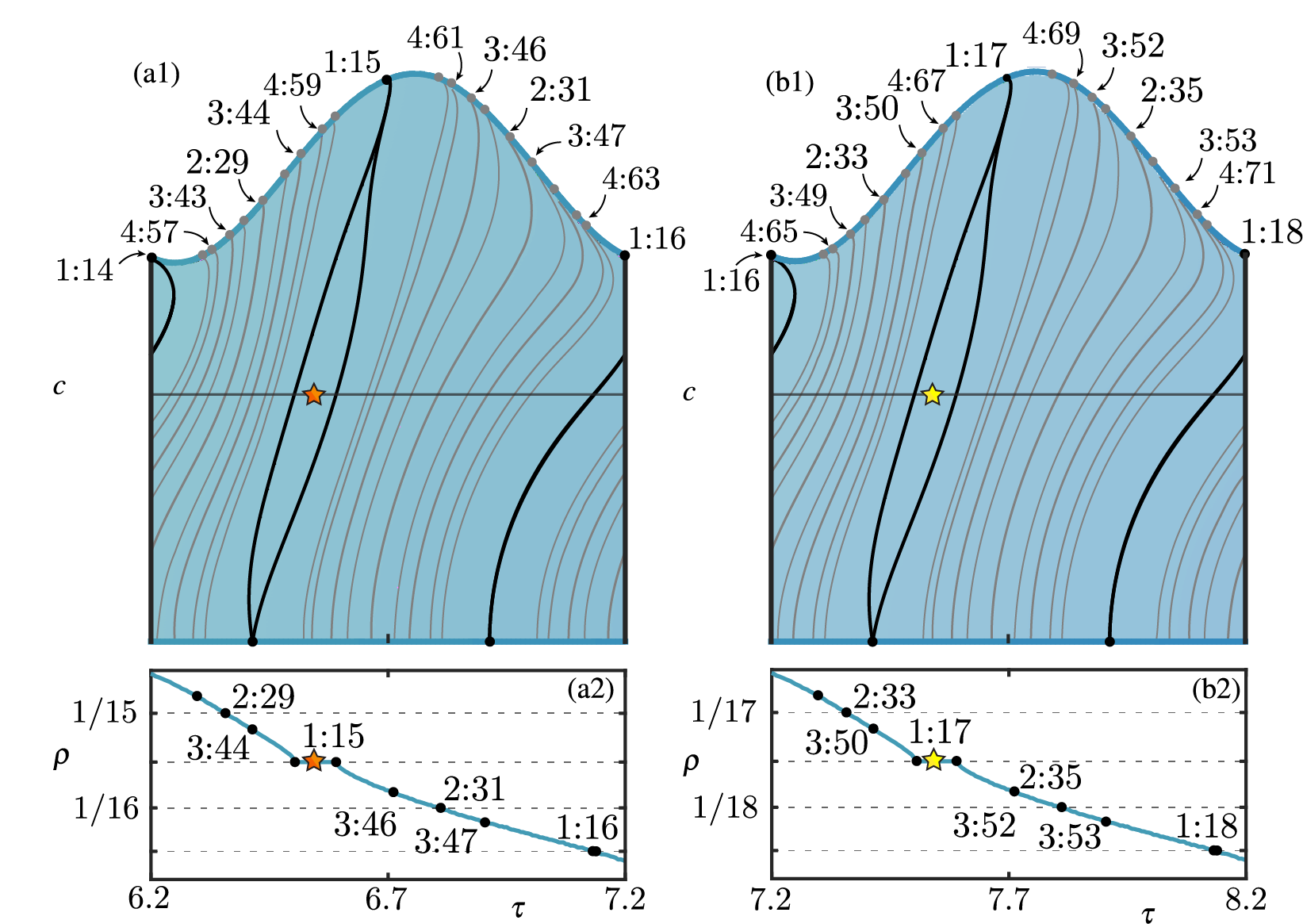}
\caption{
Preservation of the resonance structure under the shift $\tau\mapsto\tau+1$,
governed by the transformation $\Phi$, illustrated for Eq.~\eqref{eq:pss} with $\alpha=1.25$. 
Panels~(a1) and (b1) show the rotation number $\rho$ in the $(\tau,c)$--plane for $\tau\in[6.2,7.2]$ and $\tau\in[7.2,8.2]$, respectively, in the same layout as Figure~\ref{fig:tau_beta_125_2pi}. 
Panels~(a2) and (b2) show slices at fixed $c=4$, displaying the devil's staircase graph of $\rho$ where dots indicate saddle-node bifurcations. }
\label{fig:Phi_repeat}
\end{figure}

In particular, we find no 
evidence of torus breakdown near the `bifurcations' $\textsf{SNT}$ or $\textsf{GT}$: 
although theory guarantees that homoclinic tangencies and torus break-up must occur, the 
interval over which they do so is exponentially small in the frequency ratio. Owing to the 
very small rotation number in this region, we observe that those structures near \textsf{SNT} and \textsf{GT} lie beyond numerical 
detectability. The loci $\textsf{SNT}$ and $\textsf{GT}$ can therefore 
be treated as smooth curves for all practical purposes. We located 
$\textsf{SNT}$ by sweeping $\tau$ leftward for each fixed $c$ and 
identifying the rightmost value at which a torus is still found; here, we used steps of $0.0025$ in $\tau$ for 
$c$ values $0.0053$ apart. The locus $\textsf{GT}$ was traced by 
continuing a period-$50$ saddle locked orbit $\Gamma_u$ from the zero-forcing 
line through the $(\tau,c)$-plane.

\section{Reappearance of resonance tongues for large $\tau$}
\label{section:global}

\begin{figure*}[t!] 
  \centering
  \includegraphics[scale=.835]{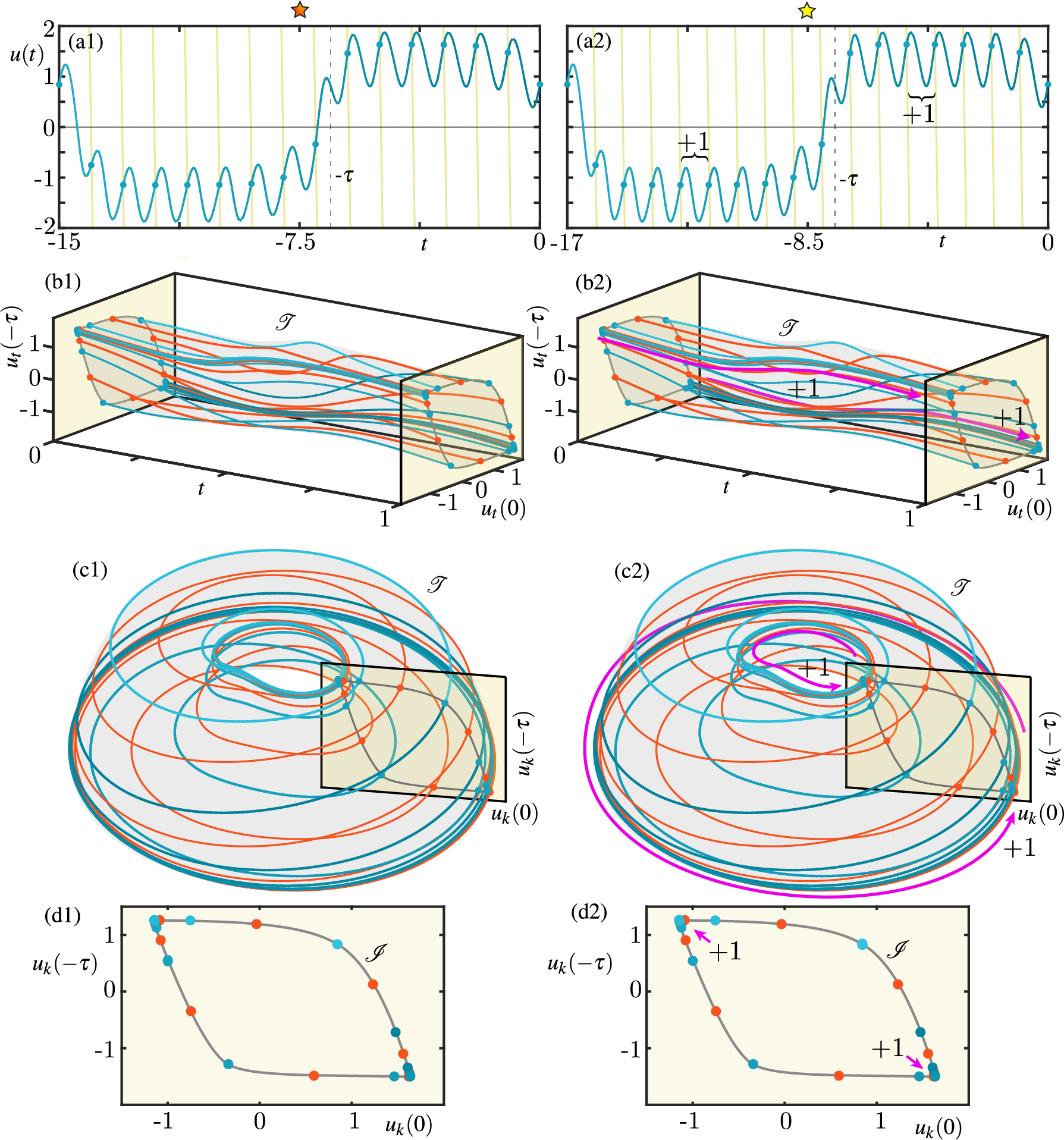}
  \caption{Locked orbits of Eq.~\eqref{eq:pss} with  $\alpha = 1.25$, $c=4$, $\omega = 2\pi$, and delay $\tau = 6.54$ in the first column and $\tau = 7.54$ in the second column. Each row illustrates representative dynamics by showing: time series over one period (a); attracting (blue) and saddle (orange) locked orbits represented in $(t, u_t(0), u_t(-\tau))$-space in between the planes $t=0$ and $t=1$ (b); a torus embedding showing these orbits, where the planes $t=0$ and $t=1$ are identified (c); and the headpoints of the locked orbit in the stroboscopic Poincar\'e section (plane at $t=0$) in panel~(d). We also show in panels~(b) and (c)  the unstable manifold of the saddle orbit, which forms the invariant torus $\mathcal{T}$ (grey), and the invariant circle $\mathcal{I}$ in panel~(d).}\label{fig:applicationdehn}
\end{figure*}

We now discuss the periodicity observed in Figures~\ref{fig:tau_beta_125_2pi} and \ref{fig:tau_beta_75_2pi}  of the resonance structure with period one in $\tau$, for sufficiently large $\tau$.
Figure~\ref{fig:Phi_repeat} shows the bifurcation set and the rotation number $\rho$ of Eq.~\eqref{eq:pss} in the $(\tau,c)$-plane with $\alpha = 1.25$ and $\omega = 2\pi$ for $\tau\in[6.2,7.2]$ in panel~(a1) and for $\tau\in[7.2,8.2]$ in panel~(b1), respectively; both are beyond $\tau$-range in Figure~\ref{fig:tau_beta_125_2pi}. We overlay all $p$:$q$ resonance tongues of the form $1$:$q$ up to $5$:$q$, together with their roots on the torus bifurcation curves $\textsf{T}$ and $\textsf{T}_0$, with the rotation number $\rho$ colored according to the colorbar provided in Figure~\ref{fig:tau_beta_75_2pi}. Figures~\ref{fig:Phi_repeat}(a2)--(b2) show associated slices at fixed $c=4$, displaying the graph of $\rho$ where the dots indicate the associated resonance tongues.

Comparing Figure~\ref{fig:Phi_repeat}(a1) and (b1), we observe a near-identical arrangement of the resonance tongues under the shift $\tau\mapsto\tau+1$, where a $p$:$q$ resonance tongue appears transformed to 
\begin{equation}
    \label{eq:Phi}
    \Phi\!\left(p{:}q\right)=p{:}(q+2p).
\end{equation}
For example, the $1{:}14$, $2{:}29$, and $3{:}44$ resonance tongues are mapped to $1{:}16$, $2{:}33$, and $3{:}50$, respectively. 
Accordingly, the graph of $\rho$ in panel~(a2) aligns with that in panel~(b2), subject to \eqref{eq:Phi}.

\smallskip

To illustrate how the transformation $\Phi$ manifests at the level of individual solutions, Figure~\ref{fig:applicationdehn} compares periodic orbits of Eq.~\eqref{eq:pss} for fixed $\alpha =1.25$, $c=4$ and $\omega = 2\pi$, for the two delays $\tau=6.54$ and $\tau=7.54$, shown in the left and right columns.
These periodic orbits lie within the $1{:}15$ and $\Phi(1{:}15)=1{:}17$ resonance tongues, respectively (indicated by the orange and yellow stars in Figure~\ref{fig:Phi_repeat}). 

Row~(a) of Figure~\ref{fig:applicationdehn} displays the time series of the attracting periodic orbit over 
one period, with vertical lines marking each forcing period. Rows~(b)--(d) follow the same format as Figure~\ref{fig:rotationumber_plot}: row~(b) shows the attracting and saddle locked orbits together with the invariant torus $\mathcal{T}$ in the $(t, u_t(0), u_t(-\tau))$-space, with the planes $t=0$ and $t=1$; row~(c) shows the torus embedded into $\mathbb R^3$ by~\eqref{eq:toruscoordinates}; and row~(d) the associated headpoints in the $(u_t(0),u_t(-\tau))$-plane with the invariant circle $\mathcal{I}$.

The time series in Figure~\ref{fig:applicationdehn}(a2) differs from that in panel~(a1) by two additional oscillations within the period. Accordingly, there are two additional segments between the planes at $t=0$ and $t=1$ in panel~(b2) compared to (b1). The locked periodic orbit in panel~(c2) winds two additional times around the longitude direction of $\mathcal{T}$ compared to the one in panel~(c1), corresponding to two additional headpoints in panel~(d2) compared to (d1).

\medskip

The key 
observation is that $\Phi$ admits a natural algebraic interpretation via 
the identification of the reduced fraction $p{:}q$ with the integer vector 
$(p,q)\in\mathbb{Z}\times\mathbb{Z}$, giving the lift 
\[
\widetilde{\Phi}\!\begin{pmatrix}
    p \\ q
\end{pmatrix} = M\!\begin{pmatrix}
    p \\ q
\end{pmatrix}\!,\,\,\text{with}
\,\,
M = \begin{pmatrix} 1 & 0 \\ 2 & 1 \end{pmatrix},
\]
which recovers the map \eqref{eq:Phi} through $(p,q)\mapsto p{:}q$. 
Importantly, $\widetilde{\Phi}$ is a linear map with $M\in \mathrm{SL}(2,\mathbb{Z})$, the group of $2\times2$ integer matrices with determinant one, and this means that $\Phi$ preserves Farey addition
\[
\frac{p}{q}\oplus\frac{\ell}{m}
=
\frac{p+\ell}{q+m},
\]
since it is is vector addition for $\widetilde{\Phi}$. 
Note that since \(M\in\mathrm{SL}(2,\mathbb Z)\), primitive integer vectors are mapped to primitive integer vectors, so reduced fractions are mapped to reduced fractions.

The matrix $M$ has a topological interpretation: 
$M = L^2$, where
\[
L=\begin{pmatrix}1&0\\1&1\end{pmatrix}\in\mathrm{SL}(2,\mathbb{Z})
\]
is the standard generator of a \emph{Dehn twist} about one of the 
fundamental cycles of the torus, acting as $(p,q)\mapsto(p,q+p)$; see, 
e.g., Chapter~3 of \cite{farb2011primer}. Hence, $M=L^2$ 
corresponds to two successive Dehn twists along the same cycle, and the 
passage $\tau\mapsto\tau+1$ can be interpreted algebraically as a \emph{double Dehn twist} on the 
torus. To visualize this, imagine slicing the torus along a longitude — 
much like cutting a bagel — twisting one side twice along the cut before reattaching it; see Figure~\ref{fig:applicationdehn}(c).

The map \(\Phi\) admits a natural extension from rational rotation numbers to all of \([0,1]\). Indeed, the rational map \eqref{eq:Phi} 
is given by the M\"obius transformation
\begin{equation}
    \label{eq:Phix}
    \Phi(\rho)=\frac{\rho}{1+2\rho},
\end{equation}
which is continuous and strictly increasing on \([0,1]\). Since the rational numbers are dense in \([0,1]\), this uniquely determines the action of \(\Phi\) on irrational rotation numbers as the continuous extension of the rational case.
Note that, because $0<\Phi(\rho)<\rho$ for all $\rho\in(0,1]$, the rotation number $\rho$ decreases as $\tau$ increases, in the sense that the associated resonances occur at smaller rotation numbers, which we can clearly observe in Figures~\ref{fig:tau_beta_125_2pi} and \ref{fig:tau_beta_75_2pi}.

\medskip

We can formulate our observations as follows.

\begin{proposition}
\label{prop:Phi}
[Repeating locked and quasiperiodic orbits]
Consider the pfS\&S model~\eqref{eq:pss} with $\omega=2\pi$. 
Assume that, along the zero-forcing line $c=0$, the model has for sufficiently
large $\tau$ a normally hyperbolic attracting periodic orbit $\Gamma_s(\tau)$ whose period satisfies
\begin{equation}
    \label{eq:Tasympt}
    T(\tau) = 2\tau + C(\alpha) + \mathcal{O}(1/\tau),
\qquad \text{as } \tau \to +\infty,
\end{equation}
where $C(\alpha)$ depends on $\alpha$.
Let
\[
\rho(\tau):=\frac{1}{T(\tau)} .
\]
Then, under the shift $\tau\mapsto \tau+1$,
\begin{equation}
    \rho(\tau+1)
    =
    \Phi(\rho(\tau))
    +
    \mathcal{O}(\rho(\tau)^3),
\end{equation}
where $\Phi$ is given by~\eqref{eq:Phix}.
Consequently, the map $\Phi$, obtained from the zero-forcing asymptotics,
organizes to leading order the reappearance of locking and quasiperiodic
dynamics for small periodic forcing $c>0$.
\end{proposition}
The proof of Proposition~\ref{prop:Phi} is found in Appendix~\ref{app:largedelayregime}. 

\smallskip

The appearance of the double Dehn twist follows directly from the asymptotic scaling $T(\tau) \sim 2\tau$ in Proposition~\ref{prop:Phi}, where the coefficient $2$ reflects the symmetric periodic orbit's asymptotic square-wave structure, which traverses two slow phases per cycle. We remark that numerical evidence for the GZT model suggests an effective scaling with $n = 4$, instead yielding a quadruple Dehn twist. More generally, if the intrinsic period satisfies an asymptotic relation of the form $T(\tau) \sim n\tau$, then the corresponding transformation of rotation numbers would formally be
\[
\rho \mapsto \frac{\rho}{1 + n\rho}.
\]
Its restriction to rational rotation numbers lifted to $\mathbb Z\times\mathbb Z$ is given by $(p,q) \mapsto L^n(p,q)$, corresponding to $n$ Dehn twists. We expect and conjecture that the resonance structure of any periodically forced DDE repeats for sufficiently large $\tau$, where the corresponding matrix in $\mathrm{SL}(2,\mathbb{Z})$ is a power of $L$.

Note that the reappearance phenomenon described in Proposition~\ref{prop:Phi} is fundamentally different from the well-known reappearance mechanism in DDEs~ \cite{Yanchuk2009,Giraldo_2023}. In that setting, if a periodic orbit with period $T$ exists at $\tau = \tau_0$, then an exact $T$-periodic orbit also exists near $\tau = \tau_0 + T$, often with different stability.

\begin{figure*}[t!] 
  \centering
    \includegraphics[scale=0.86]{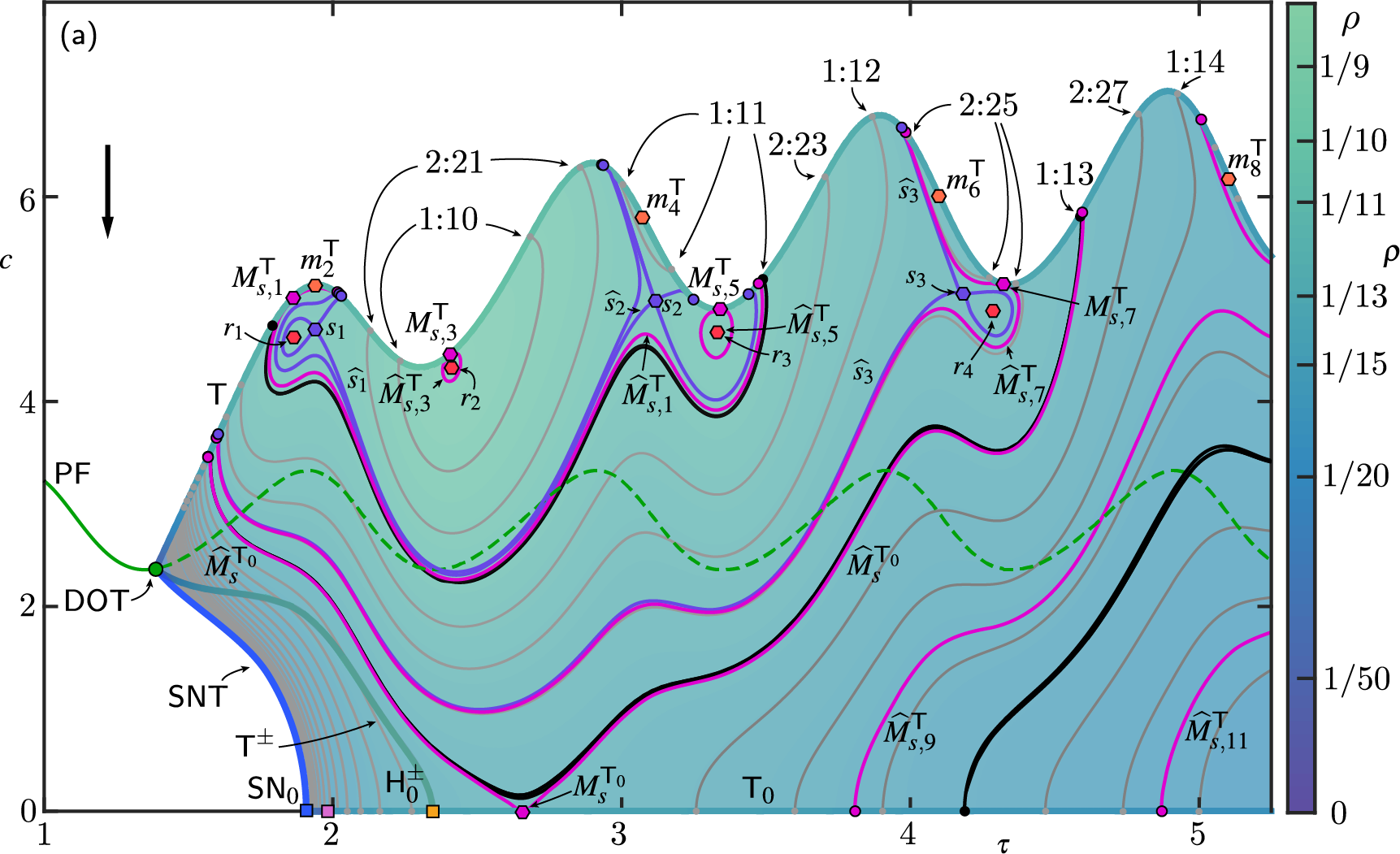}
\caption{
Bifurcation set and rotation number $\rho$ of the pfVoC model in the 
$(\tau,c)$-plane for $\alpha = 0.75$, $\omega = 2\pi$, and $\gamma = 
0.49$. The figure follows the same format as 
Figure~\ref{fig:tau_beta_75_2pi}. Along the torus bifurcation curve 
$\textsf{T}$, boundary saddles $M_{s,2k-1}^{\textsf{T}}$ (pink) and 
boundary minima $m_{2k}^{\textsf{T}}$ (orange) are marked, with the 
associated separatrices $\widehat{M}_{s,k}^{\textsf{T}}$ shown in pink. 
In the interior of the torus dynamics region, interior maxima $r_k$ 
(red) and interior saddles $s_k$ (purple) are marked, with associated 
separatrices $\widehat{s}_k$ shown in purple.
}\label{fig:pfVoC}
\end{figure*}

\section{Discussion and outlook}
\label{section:discussionandoutlook}
We carried out a comprehensive bifurcation analysis of the 
pfS\&S model~\eqref{eq:pss} with a focus on its resonance structure. 
We began by analyzing the autonomous S\&S model~\eqref{eq:ss}, providing a complete 
description of its equilibria and their stability via a linear stability analysis. With \texttt{DDE-BifTool}, we found and continued bifurcations that linear stability analysis cannot capture to complete the global bifurcation scenario. In this way, we established that, for sufficiently large delays, the sole attracting invariant structures are periodic orbits.

To assess the effect of seasonal locking, we introduced periodic forcing with angular frequency $\omega = 2\pi$ (corresponding to annual forcing, for time in years), leading to the pfS\&S  DDE model \eqref{eq:pss}. In such a periodically forced system, the intrinsic oscillations resulting from the delayed feedback interact with the external forcing, giving rise to attracting invariant tori and resonance phenomena.
We extended to the periodically forced DDE setting a method originally developed for periodically forced planar ODEs \cite{bailie2024resonancetongues,John} to compute the rotation number $\rho$ for an attracting invariant torus. A key part of this extension involves carefully handling projection issues specific to the infinite-dimensional phase space of a DDE.

We analyzed the bifurcation diagram and resonance structure of the pfS\&S model in the $(\tau,c)$-plane for 
representative values of the feedback strength: the mathematically convenient case 
$\alpha = 1.25$ and the physically more relevant case $\alpha = 0.75$. For both, resonance tongues organize around critical points of the rotation number $\rho$ --- both along the torus bifurcation curves and in the interior of 
the region with an invariant torus; this is consistent with recent observations~ \cite{bailie2024resonancetongues,John}. We found that the resonance structure repeats for large delays as described by  the rational map $\Phi$. At the level of individual 
locked orbits this map $\Phi$ encodes additional winding along the longitudinal direction of the invariant 
torus; mathematically, it is a map on rational rotation numbers that respects the Farey tree. 
For $\alpha = 1.25$, there is always a unique global attractor. However, for $\alpha = 0.75$, we additionally found bistability between 
period-one orbits and invariant tori. The associated hysteresis scenario 
involves complex bifurcation sequences, including `saddle-node' and `gluing' bifurcations of tori. Importantly, neither regime exhibits resonance overlap and, consequently, no chaotic dynamic is observed. 

\smallskip

As we mentioned in the introduction, Suarez and Schopf~ \cite{SuarezSchopf} 
derived their S\&S model in an ad hoc way, but recently Falkena \emph{et 
al.}~ \cite{falkena2019derivation} derived two S\&S-type models from a 
spatially extended ENSO model more formally. The Mori--Zwanzig model 
lacks stable periodic orbits for most physically realistic parameter 
values and is therefore unsuitable for modeling ENSO variability; we 
focus instead on the variation-of-constants (VoC) model:
\begin{equation}
\label{eq:VoC_conc}
u'(t) = u(t) - u^{3}(t) - \alpha u(t-\tau)\bigl(1 - \gamma u^{2}(t)\bigr),
\end{equation}
which reduces to the classical S\&S model~\eqref{eq:ss} when $\gamma=0$. 
For $\gamma=0.49$, the autonomous VoC model has the same qualitative 
dynamical landscape as the S\&S model, with bifurcation curves shifted 
towards larger $\tau$; see \cite{falkena2019derivation}.

A natural question is how far the findings for the pfS\&S model carry 
over to a more physically grounded model. To address this, we briefly 
examine the periodically forced VoC (pfVoC) model, obtained by adding 
the forcing term $c\cos(2\pi t)$ to~\eqref{eq:VoC_conc}. 
Figure~\ref{fig:pfVoC} shows the bifurcation diagram of the pfVoC model in the $(\tau,c)$-plane for $\gamma=0.49$, $\alpha=0.75$, and 
$\omega=2\pi$, in the format of  Figure~\ref{fig:tau_beta_75_2pi}. The pfVoC model shares the same global 
bifurcation structure as the pfS\&S model with $\alpha=0.75$: the region 
of attracting tori is bounded by the curve of torus bifurcations $\textsf{T}$, $\textsf{T}_0$, and the 
saddle-node of tori $\textsf{SNT}$; resonance tongues do not overlap so 
that $\rho$ is continuous throughout; and the reappearance mechanism 
governed by~\eqref{eq:Phix} persists for large $\tau$, though it is  less apparent in Figure~\ref{fig:pfVoC}.
This supports the conjecture that the repeating resonance structure for large $\tau$ is a typical phenomenon of periodically forced DDEs.

Despite this structural similarity, the resonance organization of the 
pfVoC model is considerably richer: the graph of $\rho=\rho(\tau,c)$ 
contains more critical points, with an alternating sequence of maxima 
and minima of $\rho$ along the torus bifurcation curve $\textsf{T}$, 
and numerous interior maxima and interior saddles in the region of 
attracting tori, producing a more intricate arrangement of resonance 
tongues and separatrices; see Figure~\ref{fig:pfVoC} and its caption 
for details.

\smallskip

A further point of discussion concerns the role of the forcing timescale 
on the dynamics of the pfS\&S model. As noted in the introduction, the timescale of the model is not fixed a priori and depends on many climatological parameters. On the other hand, varying $\omega$ while keeping the period in physical years fixed rescales the delay, the feedback strength, the forcing amplitude, and, more consequentially, the strength of the nonlinearities, which can promote resonance tongues overlap and lead to chaotic dynamics~ \cite{BolducIFAC}. This is consistent with existing results: Anikushin and Romanov~ \cite{ANIKUSHIN2023133653} reported positive Lyapunov exponents for trajectories at $\omega = 1$, and Oishi and Sekine~ \cite{Oishi2021} proved the existence of large overlapping resonance tongues for broad ranges of $\omega$. Understanding how resonance 
tongues grow, interact, and overlap as $\omega$ decreases, and how the graph of $\rho$ evolves in this process, appears to be essential for determining whether periodically forced delay models 
can reproduce realistic ENSO variability. 
In this respect, the GZT model offers a useful comparison: its timescale is fixed by construction, but it contains an explicit parameter controlling the strength of the nonlinearity, which allows one to tune the degree of resonance overlap directly. The pfS\&S model lacks such a parameter, and the role of $\omega$ as a proxy for nonlinearity strength is a modeling subtlety that deserves further investigation.

\smallskip

Other promising extensions of the model include breaking the $\mathbb Z_2$-symmetry of the S\&S framework to capture the observed asymmetry between El Niño and La Niña events~ \cite{an2004nonlinearity,larkin2002enso}, 
incorporating state-dependent delays following the approach of Keane \emph{et al.}~ \cite{KKstatedep}, 
and including stochastic forcing, which has been shown to enhance irregular variability in periodically forced delay models~ \cite{zivkovic2013enso}.

\section*{Acknowledgments}    
We thank Stefan Ruschel, Courtney Quinn and Henk A. Dijkstra for their advice and insightful discussions. 
This research was supported by Royal Society Te Ap\={a}rangi Marsden Fund grant \#19-UOA-223.

\section*{Data Availability Statement}

All data are available within the article or supplemental information.

\appendix

\section{Proof of Proposition~\ref{prop:hopf}}
\label{appendix:local}
The equilibria of \eqref{eq:ss} satisfy 
\[
p-p^3+\alpha p=0,
\]
so the trivial equilibrium $p_0=0$ exists for all $(\alpha,\tau)$, while
\[
p^\pm=\pm\sqrt{1-\alpha}
\]
exist only when $\alpha<1$.

For $\alpha=1$, the equilibria $p^\pm$ collide with $p_0$. Because the
system possesses a $\mathbb{Z}_2$-symmetry, this is a
pitchfork bifurcation $\textsf{P}(\tau,\alpha)$ along the line given by $\alpha = 1$, which shows (ii). 
The genericity conditions (involving higher-order terms) are standard
and generically satisfied; see Kuznetsov \cite{Kuznetsov}. 

\smallskip

Linearising \eqref{eq:ss} about an equilibrium $p$ gives
\[
y'(t)=(1-3p^2)\,y(t)-\alpha\,y(t-\tau).
\]

\noindent
At the trivial equilibrium $p_0$ this yields
\begin{equation}
\label{eqn:linearise}
y'(t)=y(t)-\alpha y(t-\tau),
\end{equation}
with characteristic equation
\begin{equation}
\label{eqn:characterizticequationu0}
\lambda = 1-\alpha e^{-\lambda\tau}.
\end{equation}

At the equilibria $p^\pm$ we obtain
\[
y'(t)=(3\alpha-2)y(t)-\alpha y(t-\tau),
\]
with characteriztic equation
\begin{equation}
\label{eq:charpp}
\lambda=3\alpha-2-\alpha e^{-\lambda\tau}.
\end{equation}

\smallskip

For $0<\alpha<1$, the characteriztic equation
\eqref{eqn:characterizticequationu0} has one real positive root
$\lambda_1>0$ and one real negative root $\lambda_2<0$, with all
remaining roots having real part less than $\lambda_2$, as shown in \cite{ANIKUSHIN2023133653}. Hence $p_0$ is a saddle with a
one-dimensional unstable manifold $W^u(p_0)$.

\smallskip

For $\alpha>1$, we seek purely imaginary roots
$\lambda=i\nu$ of \eqref{eqn:characterizticequationu0}. Substituting
$\lambda=i\nu$ gives
\[
\cos(\nu\tau)=\frac{1}{\alpha}, \qquad
\sin(\nu\tau)=\frac{\nu}{\alpha}.
\]
Using the Pythagorean trigonometric identity yields $\nu^2=\alpha^2-1$
which is positive for $\alpha>1$. Writing
\[
\nu_0=\sqrt{\alpha^2-1},
\]
we obtain the Hopf curves $\textsf{H}_k(\tau_k,\alpha)$ where
\begin{equation}
    \label{eq:tauk_app}
    \tau_k
=
\frac{\arccos(1/\alpha)+2\pi k}{\nu_0},
\qquad k\in\mathbb{N}_0.
\end{equation}
Each $\tau_k$ corresponds to a pair of eigenvalues
$\pm i\nu_0$ crossing the imaginary axis. The principal branch
$\tau_0$ marks the first loss of stability of $p_0$. Because $p_0$ is invariant under the $\mathbb{Z}_2$-involution, the emanating periodic orbits satisfy \eqref{eqn:symmetry}.

For $k\geq 1$, as $\alpha\to 1^+$, the numerator of \eqref{eq:tauk_app}
converges to $2\pi k$, whereas $\nu_0\to 0^+$,
so $\tau_k\to+\infty$. Alternatively, as $\alpha\to+\infty$, the numerator converges to $\pi/2+2\pi k$
and $\nu_0\to+\infty$, so $\tau_k\to 0$. This proves (i) and (iv).

\smallskip

For $\alpha<1$ we analyse the characteriztic equation
\eqref{eq:charpp}. Substituting $\lambda=i\nu$ and separating real and
imaginary parts gives
\begin{align}
\cos(\nu\tau)&=\frac{3\alpha-2}{\alpha}, \\
\sin(\nu\tau)&=\frac{\nu}{\alpha},
\end{align}
which yields
\[
\nu^2=\alpha^2-(3\alpha-2)^2
=4(1-2\alpha)(\alpha-1).
\]
This quantity is positive only for $\tfrac12<\alpha<1$. Writing
\[
\nu_0=2\sqrt{(1-2\alpha)(\alpha-1)},
\]
we obtain
\begin{equation}
    \label{eq:taukpm_app}
    \tau_k^\pm
=
\frac{\arccos\!\big(\tfrac{3\alpha-2}{\alpha}\big)+2\pi k}{\nu_0},
\qquad k\in\mathbb{N}_0.
\end{equation}
Each such value corresponds to a pair of purely imaginary eigenvalues
$\pm i\nu_0$. The principal branch $\tau_0^\pm$ marks the first loss
of stability of $p^\pm$ as $\tau$ increases. The emanating periodic orbits appear in symmetric pairs related by \eqref{eqn:symmetry}, corresponding to the two equilibria $p^\pm$.

For $k\geq 1$, as $\alpha\to 1^-$, the numerator of \eqref{eq:taukpm_app}
converges to $2\pi k$,
whereas $\nu_0\to 0^+$, so $\tau_k^\pm\to+\infty$. As $\alpha\to\tfrac{1}{2}^+$, the numerator tends to
$\pi(2k+1)$, whereas $\nu_0\to 0^+$, so $\tau_k^\pm\to+\infty$. This proves (iii) and (v).

\smallskip

Finally, rewriting \eqref{eqn:characterizticequationu0} as
\[
f(\lambda;\alpha,\tau)=1-\alpha e^{-\tau\lambda}-\lambda,
\]
one checks that at $(\alpha,\tau)=(1,1)$ the root $\lambda=0$ has
multiplicity two. Indeed,
\[
\lim_{\lambda\to0}\frac{f(\lambda;1,1)}{\lambda}=0,
\qquad
\lim_{\lambda\to0}\frac{f(\lambda;1,1)}{\lambda^2}=-\frac12,
\]
so the algebraic multiplicity is exactly two, and this establishes~(vi).

We do not verify the Hopf and pitchfork nondegeneracy conditions
analytically; these have been confirmed numerically in later sections using
\texttt{DDE-BifTool}.

\section{Proof of Proposition~\ref{prop:Phi}}\label{app:largedelayregime}

The period law~\eqref{eq:Tasympt} is motivated by a large-delay
square-wave analysis of Chow--Mallet-Paret type: after rescaling time by
$s=t/\tau$, the zero-forcing equation becomes a singularly perturbed delay
equation whose outer solutions lie near the plateaus
$\pm\sqrt{1+\alpha}$, with rapid transitions between them. This formal picture
suggests a period $T(\tau)$ of the form $2\tau+C(\alpha)$ plus some small corrections, 
where the constant $C(\alpha)$ is determined by an associated heteroclinic
matching problem. In our computations, solving the equation with
\texttt{dde23} for $\tau\in[50,500]$ and fitting the measured periods to
$T(\tau)=a_1\tau+a_2+{a_3}/{\tau}$ 
consistently gives $a_1\approx2$, with $a_2$ depending on $\alpha$ and with only
small $1/\tau$ corrections. This supports the asymptotic form 
\eqref{eq:Tasympt} which gives
\[
T(\tau+1)
=
2(\tau+1)+C(\alpha)+\mathcal{O}(1/\tau)
=
T(\tau)+2+\mathcal{O}(1/\tau),
\]
because both remainders at $\tau$ and $\tau+1$ are of order
$\mathcal{O}(1/\tau)$. Hence, 
\[
\rho(\tau+1)
=
\frac{1}{T(\tau)+2+\mathcal{O}(1/\tau)}.
\]
Since $T(\tau)=\mathcal{O}(\tau)$ and $\rho(\tau)=1/T(\tau)$, this becomes
\[
\rho(\tau+1)
=
\frac{\rho(\tau)}
     {1+2\rho(\tau)+\mathcal{O}(\rho(\tau)^2)}.
\]
Expanding the denominator gives
\[
\rho(\tau+1)
=
\frac{\rho(\tau)}{1+2\rho(\tau)}
+
\mathcal{O}(\rho(\tau)^3).
\]
Thus, to leading order,
\[
\rho(\tau+1)=\Phi(\rho(\tau)),
\quad\text{with } \Phi(\rho)=\frac{\rho}{1+2\rho}.
\]

Assuming that the periodic orbit $\Gamma_s$ is normally hyperbolic, it
persists under sufficiently small forcing $c>0$. Resonance tongues are then
expected to emanate from rational values $\rho=p/q$ along the
$c=0$ line, corresponding to $p\!:\!q$ phase locking. Between these,
one expects quasiperiodic dynamics with irrational rotation number. Hence, the
map $\Phi$, although derived from the zero-forcing asymptotics, gives the
leading-order organization of the nearby forced dynamics.

\bibliographystyle{plain}  %
   \bibliography{references}

\end{document}